# Adaptive Nonlinear Model Predictive Control of Monoclonal Antibody Glycosylation in CHO Cell Culture

Yingjie Ma[a], Jing Guo[a,b], Alexis B. Dubs[a], Krystian K. Ganko[a] and Richard D. Braatz[a,*]

[a]*Massachusetts Institute of Technology, Cambridge, MA 02139, USA*
[b]*Polytechnique Montréal, Montréal, QC H3T 0A3, Canada*

ABSTRACT

N-glycosylation is a critical quality attribute of monoclonal antibodies (mAbs), the dominant class of biopharmaceuticals. Controlling glycosylation remains difficult due to intrinsic pathway complexity, limited online measurements, and a lack of tailored control strategies. This work applies an adaptive nonlinear model predictive control (ANMPC) framework to a fed-batch mAb production process, using a multiscale model that links extracellular conditions to intracellular Golgi reactions to predict glycan profiles. Model parameters are updated online as new measurements arrive, after which a shrinking-horizon optimization computes the control inputs; only the first control move is implemented each cycle. Case studies show that, with a minimal day-1 galactose excitation, ANMPC mitigates model–plant mismatch and achieves up to 130% and 96% higher performance than open-loop optimization and state NMPC, respectively. Under more realistic conditions (partial measurement availability and longer preparation time), ANMPC maintains comparable performance, indicating robustness to practical limitations. Overall, the results demonstrate that ANMPC can actively shape glycan distributions in silico and offers a viable path toward closed-loop control of mAb glycosylation.

## 1. Introduction

Monoclonal antibodies (mAbs) are among the highest-revenue biopharmaceuticals, representing US \$217 billion of the global biopharmaceutical market (US \$343 billion) in 2021 [46]. They are widely administered to treat cancer, infectious diseases, and inflammatory disorders [31]. The N-linked glycans attached to the Fc region of a mAb profoundly influence its bioactivity and therapeutic efficacy [16, 27]. Because the glycan distribution must remain within tight limits, regulatory guidelines classify glycosylation as a critical quality attribute (CQA) for mAb products [2]. Variations in core fucosylation and terminal galactosylation, for example, directly modulate antibody-dependent cellular cytotoxicity (ADCC) and complement-dependent cytotoxicity (CDC)—the primary mechanisms through which mAbs destroy target cells [15, 29, 38]. Conversely, an increased proportion of high-mannose glycans accelerates serum clearance and reduces therapeutic potency [9].

Although essential for mAb bioactivity and therapeutic efficacy, achieving a consistent batch-to-batch glycan profile remains a formidable challenge for biopharmaceutical manufacturers [7, 30, 40]. Chinese hamster ovary (CHO) cells are the industry workhorse for mAb production because they perform complex post-translational modifications and secrete proteins that are both human-compatible and biologically active [21]. The CHO-based production workflow comprises a sequence of template-driven steps—DNA transcription, mRNA translation—followed by protein folding and glycosylation, which, unlike the earlier stages, do not rely on a nucleic acid template [39, 41]. Instead, glycosylation emerges from a non-template enzymatic network of thousands of coupled reactions in the endoplasmic reticulum and Golgi apparatus, all modulated by diverse intracellular factors [19]. Consequently, glycan biosynthesis and conjugation to the antibody frequently yield a heterogeneous mixture of glycoforms [41]. Overall, the inherent complexity of the glycosylation process, the lack of robust real-time glycan analytics, and the absence of dedicated control strategies have so far limited the implementation of online glycosylation control [41].

Most published control strategies for mAb bioprocesses have targeted the productivity of mAb or the concentrations of extracellular metabolites (e.g., glucose, glutamine, lactate) rather than the glycan profile itself [4, 5, 23, 35, 36]. The emphasis on productivity reflects the historical absence of mechanistic models that couple extracellular conditions to the intracellular glycosylation network—yet maximizing titer (mAb concentration) can be at odds with meeting a specified glycan distribution. Seeking to reconcile this trade-off, Ref. [47] investigated how dynamic feed profiles influence N-glycan quality, but their policies were still framed in terms of holding nutrient concentrations at preset levels rather than directly steering the glycosylation state space. Reference [48] introduced a model predictive control (MPC) scheme to maintain a target high-mannose (HM) fraction, but the underlying model only comprised a single manipulated variable and a single glycoform, limiting its ability to balance yield enhancement with multi-attribute quality specifications.

The emergence of multiscale mechanistic models that explicitly link extracellular metabolism to intracellular glycosylation has enabled high-fidelity open-loop optimizations capable of simultaneously boosting product titer and shaping the glycan profile [17, 21]. Reference [21] validated their optimization experimentally, achieving an over 90% increase in the desired glycoform and underscoring the promise of this approach. Nevertheless, both studies generated feed trajectories offline; such open-loop policies are vulnerable to model–plant mismatch and unanticipated disturbances [34]. Embedding the resulting partial-differential–algebraic equation (PDAE) models directly in a nonlinear MPC (NMPC) framework remains computationally prohibitive—a single optimization may require 20 h [26] to 40 h [17]. To reduce this burden, several groups have linearized the multiscale model into a process-gain matrix and applied the resulting controller in fed-batch [24] and perfusion bioreactors [28], but such linearizations are accurate only near the operating point used for their construction [41], limiting their robustness for the wide excursions typical of industrial fed-batch or disturbance-prone perfusion processes.

*Corresponding author
yingjma@mit.edu (Y. Ma); jing.guo@polymtl.ca (J. Guo); adubs@mit.edu (A.B. Dubs); kkganko@mit.edu (K.K. Ganko); braatz@mit.edu (R.D. Braatz)
ORCID(s): 0009-0006-4458-9859 (Y. Ma); 0000-0003-4304-3484 (R.D. Braatz)

A second challenge is parameter identification for the multiscale model. Adaptive MPC techniques—updating model parameters online from process data—have been used to address model uncertainty with promising results [11, 14, 33]. Motivated by these gaps, we develop an adaptive NMPC (ANMPC) controller for a fed-batch mAb bioreactor using the multiscale PDAE glycosylation model of Ref. [21]. At each sampling time, ANMPC updates model parameters using all available measurements and then recomputes the control moves via model-based optimization. To solve the resulting dynamic optimization and parameter estimation problems efficiently and robustly, we employ control-vector parameterization (CVP) with embedded simulations accelerated by a parallel quasi-steady-state (QSS) approach proposed in our previous work [26]. Case studies initialized from multiple parameter sets compare open-loop optimization, state NMPC, and ANMPC; ANMPC mitigates model–plant mismatch and achieves up to 130% higher penalized merit than the alternatives. We further evaluate ANMPC under more realistic measurement availability (no NSD measurements) and with a 4 h analytics/actuation delay, observing no material performance loss, which demonstrates robustness.

In what follows, Section 2 describes the multiscale glycosylation model; Section 3 formulates the control problem; Section 4 presents the ANMPC framework; Section 5 examines actuator choices and demonstrates the algorithm through case studies; and Section 6 concludes.

## 2. Glycosylation model

Figure 1 illustrates the three-level multiscale glycosylation model employed in this work:

(1) Bioreactor-level cell culture model – predicts viable cell density, extracellular metabolite concentrations, specific productivity, and the concentrations of secreted mAb.

(2) Intracellular nucleotide sugar donor (NSD) synthesis model – converts the extracellular metabolite information from the culture model into intracellular NSD using a detailed synthesis pathway and associated rate laws.

(3) Golgi-level glycosylation reaction model – receives the NSDs as sugar donors and computes the fractional distribution of glycoforms leaving the secretory pathway via a network of enzyme-catalyzed reactions coupled with glycoprotein transport equations.

The overall structure follows Ref. [21] and is extended to incorporate manganese- and ammonia-dependent kinetics reported by Ref. [45]. Information propagates downstream from the extracellular environment through the NSD synthesis layer and to the Golgi reaction network, and the resulting NSD consumption and glycoform flux then feeds back to the NSD synthesis model and cell culture model—closing the loop and capturing the two-way coupling between reactor conditions and intracellular glycosylation.

The multiscale model is formulated as a large-scale PDAE system that couples 30 ordinary differential equations (ODEs), 34 partial differential equations (PDEs), and numerous highly nonlinear algebraic relations. The remainder of this section introduces each submodel in turn.

### 2.1. Cell culture model

At the reactor scale, we employ an unstructured differential algebraic equation (DAE) model to describe cell growth, death, and metabolism [18, 21]. The dynamic balances for cell culture volume and cell populations are

$$\frac{dV}{dt} = F_{\text{in}} - F_{\text{out}} \tag{1}$$

$$\frac{d(VX)}{dt} = \mu V X_{\text{v}} - F_{\text{out}} X \tag{2}$$

$$\frac{d(VX_{\text{v}})}{dt} = \left(\mu - \mu_{\text{death}}\right) V X_{\text{v}} - F_{\text{out}} X_{\text{v}} \tag{3}$$

where $V$ (L) is the cell culture volume, $F_{\text{in}}$ and $F_{\text{out}}$ ($\text{L}\,\text{h}^{-1}$) are the inlet and outlet flow rates, $X$ and $X_{\text{v}}$ ($\text{cells}\,\text{L}^{-1}$) denote the total and viable cell densities, and $\mu$ and $\mu_{\text{death}}$ ($\text{h}^{-1}$) are the specific cell growth and death rates, respectively.

The viability is defined as the fraction of living cells:

$$\text{viability} = \frac{X_{\text{v}}}{X} \tag{4}$$

The specific growth ($\mu$) and death ($\mu_{\text{death}}$) rates are determined by nutrient availability and by the buildup of inhibitory byproducts. Glucose and asparagine supply is growth-limiting, whereas lactate, ammonia, and uridine exert inhibitory effects—with ammonia and uridine also driving cell death. The rate expressions are

$$\mu = \mu_{\max} f_{\text{lim}} f_{\text{inh}} \tag{5}$$

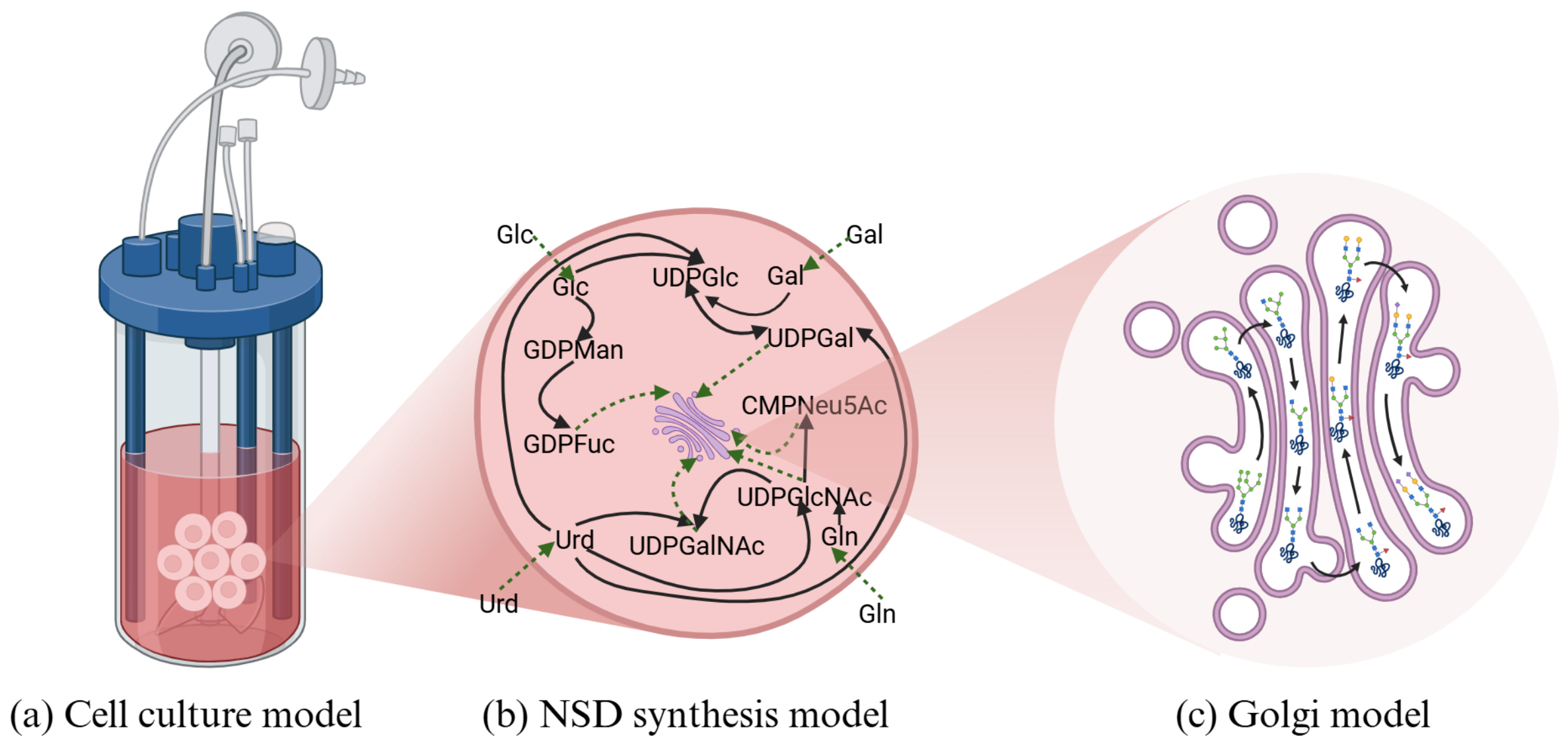

**Figure 1:** Multiscale glycosylation model

$$\mu_{\text{death}} = \mu_{\text{death,max}} \left( \frac{[\text{Amm}]}{[\text{Amm}] + K_{\text{d,Amm}}} + \frac{[\text{Urd}]}{[\text{Urd}] + K_{\text{d,Urd}}} \right) \tag{6}$$

$$f_{\text{lim}} = \frac{[\text{Glc}]}{[\text{Glc}] + K_{\text{Glc}}} \frac{[\text{Asn}]}{[\text{Asn}] + K_{\text{Asn}}} \tag{7}$$

$$f_{\text{inh}} = \frac{\text{KI}_{\text{Amm}}}{[\text{Amm}] + \text{KI}_{\text{Amm}}} \frac{\text{KI}_{\text{Lac}}}{[\text{Lac}] + \text{KI}_{\text{Lac}}} \frac{\text{KI}_{\text{Urd}}}{[\text{Urd}] + \text{KI}_{\text{Urd}}} \tag{8}$$

where $\mu_{\text{max}}$ and $\mu_{\text{death, max}}$ ($\text{h}^{-1}$) are the maximum specific growth and death rates, respectively; $f_{\text{lim}}$ and $f_{\text{inh}}$ are the substrate-limiting and metabolite-inhibiting factors, respectively; $K_{\text{Glc}}$ and $K_{\text{Asn}}$ ($\text{mmol L}^{-1}$) are the Monod half-saturation constants for glucose and asparagine; $K_{\text{d,Amm}}$ and $K_{\text{d,Urd}}$ ($\text{mmol L}^{-1}$) are the death constants for ammonia and uridine; $K_{\text{I,Amm}}$, $K_{\text{I,Lac}}$, and $K_{\text{I,Urd}}$ ($\text{mmol L}^{-1}$) are inhibition constants for ammonia, lactate, and uridine; and [Glc], [Asn], [Amm], [Lac], and [Urd] ($\text{mmol L}^{-1}$) denote the extracellular concentrations of glucose, asparagine, ammonia, lactose and uridine, respectively.

The extracellular mass balance for each metabolite is described by

$$\frac{d(V[\text{Met}])}{dt} = F_{\text{in}}[\text{Met}]_{\text{in}} - F_{\text{out}}[\text{Met}] + q_{\text{Met}} V X_{\text{v}} \tag{9}$$

where $[\text{Met}]_{\text{in}}$ is the metabolite concentration in the feed stream, [Met] is its concentration in the culture, and their units are $\text{mg L}^{-1}$ for mAb and $\text{mmol L}^{-1}$ for the other metabolites. $q_{\text{Met}}$ is the cell-specific production or consumption rate, which is expressed in $\text{pg cell}^{-1}\,\text{h}^{-1}$ for mAb and in $\text{mmol cell}^{-1}\,\text{h}^{-1}$ for the remaining metabolites. Metabolites considered in the model include ammonia (Amm), asparagine (Asn), aspartate (Asp), glucose (Glc), galactose (Gal), glutamine (Gln), glutamate (Glu), lactose (Lac), and uridine (Urd). Figure 2 depicts the corresponding metabolic network, and the individual reaction rates are

$$q_{\text{Glc}} = \left( -\frac{\mu}{Y_{X_{\text{Glc}}}} - m_{\text{Glc}} \right) \left( \frac{K_{C_{\text{Gal}}}}{K_{C_{\text{Gal}}} + [\text{Gal}]} \right)^{n_{\text{Gal}}} \tag{10}$$

$$n_{\text{Gal}} = 1 - f_{\text{Gal}} \frac{q_{\text{Gal}}}{q_{\text{Glc}}} \tag{11}$$

$$q_{\text{Gln}} = \frac{\mu}{Y_{X_{\text{Gln}}}} + q_{\text{Amm}} Y_{\text{Gln/Amm}} \tag{12}$$

$$q_{\text{Lac}} = \left( \frac{\mu}{Y_{X_{\text{Lac}}}} - Y_{\text{Lac/Glc}}\, q_{\text{Glc}} \right) \frac{\text{Lac}_{\text{max1}} - [\text{Lac}]}{\text{Lac}_{\text{max1}}} + m_{\text{Lac}} \frac{\text{Lac}_{\text{max2}} - [\text{Lac}]}{\text{Lac}_{\text{max2}}} \tag{13}$$

$$q_{\text{Amm}} = \frac{\mu}{Y_{X_{\text{Amm}}}} - Y_{\text{Amm/Urd}}\, q_{\text{Urd}} \tag{14}$$

$$q_{\text{Glu}} = -\frac{\mu}{Y_{X_{\text{Glu}}}} \tag{15}$$

$$q_{\text{Gal}} = -\frac{\mu}{Y_{X_{\text{Gal}}}} \frac{[\text{Gal}]}{[\text{Gal}] + K_{\text{Gal}}} \tag{16}$$

$$q_{\text{Urd}} = \frac{\mu}{Y_{X_{\text{Urd}}}} \frac{[\text{Urd}]}{[\text{Urd}] + K_{\text{Urd}}} \tag{17}$$

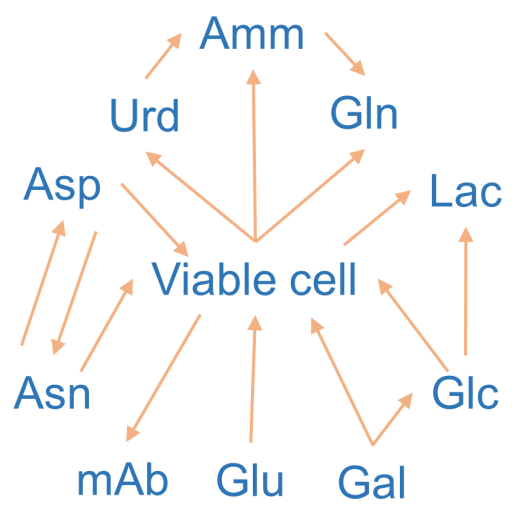

**Figure 2:** The metabolic reaction network in the cell culture. The start of an arrow refers to the reactant, while the end of the arrow denotes the product.

$$q_{\text{Asn}} = -\frac{\mu}{Y_{X_{\text{Asn}}}} - Y_{\text{Asn/Asp}}\, q_{\text{Asp}} \tag{18}$$

$$q_{\text{Asp}} = -\frac{\mu}{Y_{X_{\text{Asp}}}} - Y_{\text{Asp/Asn}}\, q_{\text{Asn}} \tag{19}$$

$$q_{\text{mAb}} = Y_{\text{mAb}_X}\mu + m_{\text{mAb}} \tag{20}$$

In these expressions, the model parameters to be estimated include the yield coefficients $Y_{X_{\text{Met}}}$ $\left(\text{cell}\,\text{mmol}^{-1}\right)$, maintenance terms $m_{\text{Met}}$ $\left(\text{mmol}\,\text{cell}^{-1}\,\text{h}^{-1}\right)$ for metabolites excluding mAb, maintenance term $m_{\text{mAb}}$ $\left(\text{pg}\,\text{cell}^{-1}\,\text{h}^{-1}\right)$ for mAb, the fraction $f_{\text{Gal}}$ (–), cross-yield ratios $Y_{\text{Met1/Met2}}$ $\left(\text{mmol}\,\text{mmol}^{-1}\right)$, saturation or inhibition constants $K_{\text{Met}}$ $(\text{mmol}\,\text{L}^{-1})$, the empirical limits $\text{Lac}_{\text{max1}}$ and $\text{Lac}_{\text{max2}}$ $(\text{mmol}\,\text{L}^{-1})$, and the product formation yield $Y_{\text{mAb}_X}$ $\left(\text{pg}\,\text{cell}^{-1}\right)$.

The extracellular mass balance equation for each glycoform $\text{GLY}_i$ is

$$\frac{d(V[\text{GLY}_i^{\text{extra}}])}{dt} = -F_{\text{out}}[\text{GLY}_i^{\text{extra}}] + V q_{\text{mAb}} X_{\text{v}} Y_i^{\text{intra}} \tag{21}$$

where $[\text{GLY}_i^{\text{extra}}]$ $\left(\text{mg}\,\text{L}^{-1}\right)$ is the extracellular concentration of glycoform $\text{GLY}_i$, and $Y_i^{\text{intra}} = [\text{GLY}_i^{\text{intra}}]/[\text{mAb}^{\text{intra}}]$ is its intracellular fractional abundance. Here, the intracellular glycoform concentration $[\text{GLY}_i^{\text{intra}}]$ $\left(\mu\text{mol}\,\text{L}^{-1}\right)$ is obtained from the Golgi model, while the total intracellular mAb level $[\text{mAb}^{\text{intra}}]$ $\left(\mu\text{mol}\,\text{L}^{-1}\right)$ is approximately $94\,\mu\text{mol}\,\text{L}^{-1}$ [20]. By analogy, the extracellular fractional abundance is $Y_i^{\text{extra}} = [\text{GLY}_i^{\text{extra}}]/[\text{mAb}]$.

## 2.2. NSD synthesis model

NSDs provide monosaccharides for the glycosylation reactions, and they are synthesized inside the cells. Therefore, the NSD synthesis model consisting of a DAE system is used to track the intracellular concentrations of seven monosaccharide donors: guanosine diphosphate mannose (GDP-Man), guanosine diphosphate fucose (GDP-Fuc), uridine diphosphate galactose (UDP-Gal), uridine diphosphate glucose (UDP-Glc), uridine diphosphate $N$-acetylgalactosamine (UDP-GalNAc), uridine diphosphate $N$-acetylglucosamine (UDP-GlcNAc), and cytidine monophosphate $N$-acetylneuraminic acid (CMP-Neu5Ac). The reaction network, depicted in Figure 1b, comprises $N_{\text{R1}}$ enzyme-catalyzed steps whose stoichiometry is encoded in the matrix $\nu^{\text{nsd}}$. For each donor $i$ $(i = 1, \dots, N_{\text{NSD}} = 7)$, the mass balance is

$$\frac{d([\text{NSD}_i^{\text{intra}}])}{dt} = \sum_{j=1}^{N_{R1}} \nu_{i,j}^{\text{nsd}} r_j^{\text{nsd}} - f_{\text{NSD}_i}^{\text{hcp/lipid}} - f_{\text{NSD}_i}^{\text{precursor}} - f_{\text{NSD}_i}^{\text{glyc}}, \quad i = 1, 2, \dots, N_{\text{NSD}} \tag{22}$$

where $[\text{NSD}_i^{\text{intra}}]$ $(\text{mmol}\,\text{L}^{-1})$ is the intracellular concentration of donor $i$; $r_j^{\text{nsd}}$ $\left(\text{mmol}\,\text{L}^{-1}\,\text{h}^{-1}\right)$ is the rate of a reaction $j$; and the sink terms $f_{\text{NSD}_i}^{\text{hcp/lipid}}$, $f_{\text{NSD}_i}^{\text{precursor}}$, and $f_{\text{NSD}_i}^{\text{glyc}}$ $\left(\text{mmol}\,\text{L}^{-1}\,\text{h}^{-1}\right)$ are consumption rates of the donor for host cell protein and glycolipid synthesis, precursor oligosaccharide assembly, and N-linked glycosylation, respectively. All reaction rates follow Michaelis–Menten kinetics that depend on the intracellular NSD concentrations together with extracellular glucose, galactose, and uridine concentrations and the intracellular glutamine level. The explicit rates are

$$r_1 = V_{\text{max},1}\frac{[\text{Gln}_{\text{intra}}]}{K_{\text{M1}_{\text{Gln}}} + [\text{Gln}_{\text{intra}}]} \tag{23}$$

$$r_{1_{\text{sink}}} = V_{\text{max},1_{\text{sink}}} \times \frac{[\text{UDPGlcNAc}]}{\left(K_{\text{M1}_{\text{sink}}} + [\text{UDPGlcNAc}]\right)\left(1 + \frac{[\text{CMPNeu5Ac}]}{\text{KI}_{1_{\text{sink}}}}\right)} \tag{24}$$

$$r_2 = V_{\text{max},2}\frac{[\text{Glc}]}{K_{\text{M2}_{\text{Glc}}} + [\text{Glc}]} \tag{25}$$

$$r_{2\text{b}} = V_{\text{max},2\text{b}}\frac{[\text{UDPGal}]}{K_{\text{M2b}_{\text{UDPGal}}}\left(1 + \frac{[\text{UDPGlcNAc}]}{K_{\text{I2A}}} + \frac{[\text{UDPGalNAc}]}{K_{\text{I2B}}} \cdots + \frac{[\text{UDPGlc}]}{K_{\text{I2C}}} + \frac{[\text{UDPGal}]}{K_{\text{I2D}}}\right) + [\text{UDPGal}]} \tag{26}$$

$$r_3 = V_{\text{max},3}\frac{[\text{Glc}]}{K_{\text{M3}_{\text{Glc}}} + [\text{Glc}]} \tag{27}$$

$$r_4 = V_{\text{max},4}\frac{[\text{UDPGlcNAc}]}{K_{\text{M4}_{\text{UDPGlcNAc}}} + [\text{UDPGlcNAc}]} \tag{28}$$

$$r_5 = V_{\text{max},5}\frac{[\text{UDPGlcNAc}]}{K_{\text{M5}_{\text{UDPGlcNAc}}}\left(1 + \frac{[\text{CMPNeu5Ac}]}{\text{KI}_5}\right) + [\text{UDPGlcNAc}]}$$

(29)

$$r_6 = V_{\max,6}\frac{[\text{UDPGlc}]}{K_{\text{M6}_{\text{UDPGlc}}}\left(1+\frac{[\text{UDPGlcNAc}]}{\text{KI}_{6\text{A}}}+\frac{[\text{UDPGalNAc}]}{\text{KI}_{6\text{B}}}\cdots+\frac{[\text{UDPGal}]}{\text{KI}_{6\text{C}}}\right)+[\text{UDPGlc}]} \tag{30}$$

$$r_{6_{\text{sink}}} = V_{\max,6_{\text{sink}}}\frac{[\text{UDPGal}]}{K_{\text{M6}_{\text{sink}}}\left(1+\frac{[\text{UDPGlc}]}{\text{KI}_{6_{\text{sink}}}}\right)+[\text{UDPGal}]}\cdots\frac{[\text{Gal}]}{[\text{Gal}]+K_{\text{regulator}}} \tag{31}$$

$$r_7 = V_{\max,7}\frac{[\text{GDPMan}]}{\left(K_{\text{M7}_{\text{GDPMan}}}+[\text{GDPMan}]\right)\left(1+\frac{[\text{GDPFuc}]}{\text{KI}_7}\right)} \tag{32}$$

$$r_{7_{\text{sink}}} = V_{\max,7_{\text{sink}}}\frac{[\text{GDPFuc}]}{K_{\text{M7}_{\text{sink}}}+[\text{GDPFuc}]} \tag{33}$$

where $V_{\max,j}$ $\left(\text{mmol L}^{-1}\,\text{h}^{-1}\right)$ is the maximum turnover rate of reaction $j$; $K_{Mj_{\text{comp}}}$ (mmol L$^{-1}$) is the Michaelis constant of a component (metabolite or NSD) in reaction $j$; $\text{KI}_j$ (mmol L$^{-1}$) is the corresponding inhibition constant; and $K_{\text{regulator}}$ (mmol L$^{-1}$) is introduced so that the factor $\frac{[\text{Gal}]}{[\text{Gal}]+K_{\text{regulator}}}$ becomes zero when extracellular Gal is absent.

The intracellular glutamine concentration (mmol L$^{-1}$) in (23) is estimated through

$$\text{Gln}_{\text{intra}} = f_{\text{Gln}}\,[\text{Gln}],$$

where $f_{\text{Gln}}$ is a constant parameter to connect intracellular Gln concentration with its extracellular counterpart.

The NSD synthesis reaction rates that depend on Urd or Gal are

$$r_{j_{\text{Urd}}} = V_{\max,j_{\text{Urd}}}\frac{[\text{Urd}]}{K_{Mj_{\text{Urd}}}+[\text{Urd}]},\quad j\in\{1,2,4,6\}, \tag{34}$$

$$r_{6_{\text{Gal}}} = \frac{V_{\max,6_{\text{Gal}}}[\text{Gal}]}{K_{M6_{\text{Gal}}}\left(1+\frac{[\text{UDPGal}]}{\text{KI}_{6\text{D}}}+\frac{[\text{Gal}]}{\text{KI}_{6\text{E}}}+\frac{[\text{Urd}]}{\text{KI}_{6\text{F}}}\right)+[\text{Gal}]}. \tag{35}$$

Additional intracellular fluxes for each $\text{NSD}_i$ are

$$f_i^{\text{hcp/lipid}} = \frac{[\text{NSD}_i]}{K_{\text{TP},i}+[\text{NSD}_i]}\frac{v_i^{\text{hcp/lipid}}\mu}{V_{\text{cell}}}, \tag{36}$$

$$f_i^{\text{precursor}} = \frac{[\text{NSD}_i]}{K_{\text{TP},i}+[\text{NSD}_i]}\frac{v_i^{\text{precursor}}q_{\text{mAb}}}{V_{\text{cell}}}, \tag{37}$$

$$f_i^{\text{glyc}} = \frac{[\text{NSD}_i]}{K_{\text{TP},i}+[\text{NSD}_i]}r_i^{\text{glyc}}, \tag{38}$$

where $K_{\text{TP},i}$ (mmol L$^{-1}$) is the transport protein saturation constant for $\text{NSD}_i$; $v_i^{\text{hcp/lipid}}$ $\left(\text{mmol cell}^{-1}\right)$ and $v_i^{\text{precursor}}$ (mmol pg$^{-1}$) are the stoichiometric requirements for host cell protein/glycolipid synthesis and precursor oligosaccharide assembly, respectively; $V_{\text{cell}}$ $\left(\text{L cell}^{-1}\right)$ is the specific cell volume; and the $\text{NSD}_i$ consumption rate in glycosylation reactions $r_i^{\text{glyc}}$ $\left(\text{mmol L}^{-1}\,\text{h}^{-1}\right)$ is [21]

$$r_i^{\text{glyc}} = \text{vel}_{\text{golgi}}\left(\frac{V_{\text{golgi}}}{V_{\text{cell}}}\right)\sum_{j=1}^{N_{\text{OS}}}\left\{v_{i,j}\left[OS_j\right](z=1)\right\} \tag{39}$$

$$\text{vel}_{\text{golgi}} = \frac{2\,q_{\text{mAb}}\times 10^{-6}}{60\,\text{MW}_{\text{mAb}}V_{\text{golgi}}\left[\text{OS}_1\right](z=0)} \tag{40}$$

where $\text{vel}_{\text{golgi}}$ $\left(\text{Golgi length min}^{-1}\right)$ is the length-normalized transit velocity, $V_{\text{golgi}}$ $\left(\text{L cell}^{-1}\right)$ is the Golgi volume, $\upsilon_{i,j}$ is the number of $\text{NSD}_i$ molecules required for one oligosaccharide $j$ molecule, $[\text{OS}_j](z)$ (mmol L$^{-1}$) is the oligosaccharide concentration at axial position $z$ (with $z=0$ at entry and $z=1$ at exit), and $\text{MW}_{\text{mAb}}$ is the molecular weight of mAb $\left(165\times10^3\,\text{g mol}^{-1}\right)$.

## 2.3. Golgi model

The Golgi model predicts the intracellular glycoform distribution using a dynamic plug flow reactor (PFR) representation of the Golgi apparatus, formulated as a PDAE system [42, 21]. The spatiotemporal balance for each oligosaccharide $\text{OS}_i$ is

$$\frac{\partial[\text{OS}_i]}{\partial t} = -\text{vel}_{\text{golgi}}\frac{\partial[\text{OS}_i]}{\partial z}+\sum_{j=1}^{N_{R2}}v_{i,j}r_j,\ i=1,\cdots,N_{\text{OS}} \tag{41}$$

where $\left[\text{OS}_i\right]$ $\left(\mu\text{mol L}^{-1}\right)$ is the local concentration of $\text{OS}_i$ in the Golgi apparatus, $r_j$ $\left(\mu\text{mol L}^{-1}\,\text{min}^{-1}\right)$ is the rate of reaction $j$, $v_{i,j}$ is the stoichiometric coefficient of $\text{OS}_i$ in reaction $j$, and $N_{R2}$ is the total number of glycosylation reactions. The glycosylation reaction network involving 33 oligosaccharides and 43 reactions is depicted in Figure 3 [45].

The glycosylation reaction rates $r_j$ fall into three mechanistic classes [42, 45]:

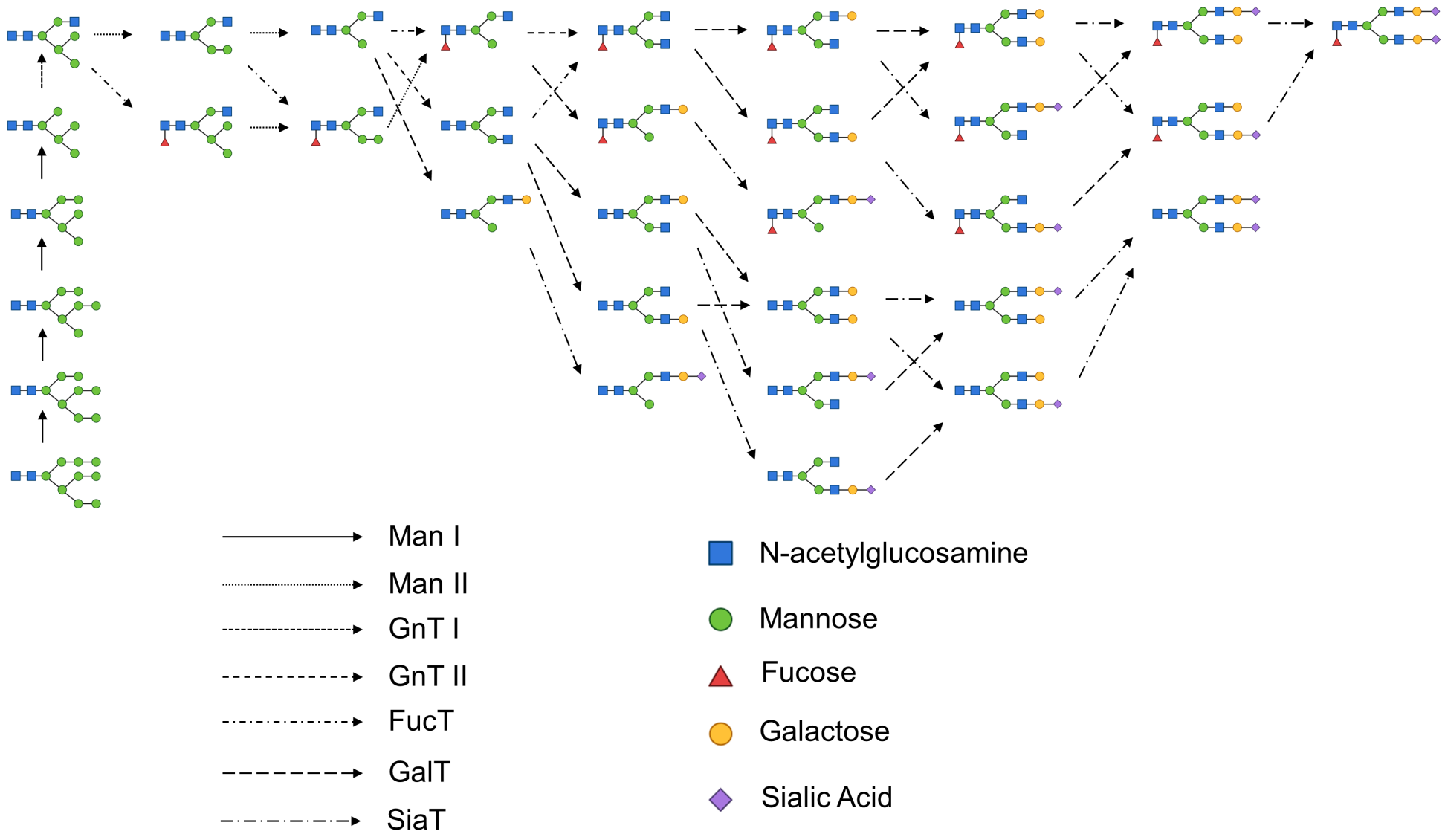

**Figure 3:** The glycosylation reaction network in the Golgi apparatus.

(i) Michaelis–Menten kinetics — used for enzymes Man I and Man II-catalyzed reactions

$$r_j = \frac{k_{f,j}\,[E_j]\,[\mathrm{OS}_i]}{K_{d,i}^{\mathrm{enz}}\left(1+\sum_{k=1}^{NC}\frac{[\mathrm{OS}_k]}{K_{d,k}^{\mathrm{enz}}}\right)} \tag{42}$$

(ii) Sequential order Bi–Bi kinetics — used for enzymes GnT I, GnT II, and GalT-catalyzed reactions

$$r_j = \frac{k_{f,j}\,[E_j]\,[\mathrm{Mn}]\,[\mathrm{NSD}_z]\,[\mathrm{OS}_i]}{K_{d,\mathrm{Mn}}^{\mathrm{enz}}\,K_{d,z}^{\mathrm{enz}}\,K_{d,i}^{\mathrm{enz}}\left(1+\frac{[\mathrm{Mn}]}{K_{d,\mathrm{Mn}}^{\mathrm{enz}}}+\frac{[\mathrm{Mn}]}{K_{d,\mathrm{Mn}}^{\mathrm{enz}}}\frac{[\mathrm{NSD}_z]}{K_{d,z}^{\mathrm{enz}}}+\frac{[\mathrm{Mn}]}{K_{d,\mathrm{Mn}}^{\mathrm{enz}}}\frac{[\mathrm{NSD}_z]}{K_{d,z}^{\mathrm{enz}}}\sum_{\substack{k=1\\k\neq i+1}}^{N_\mathrm{C}}\frac{[\mathrm{OS}_k]}{K_{d,k}^{\mathrm{enz}}}+\frac{[\mathrm{Nuc}_n]}{K_{d,n}^{\mathrm{enz}}}\frac{[\mathrm{OS}_{i+1}]}{K_{d,i+1}^{\mathrm{enz}}}+\frac{[\mathrm{Nuc}_n]}{K_{d,n}^{\mathrm{enz}}}\right)} \tag{43}$$

(iii) Random-order Bi–Bi kinetics — used for enzymes FucT- and SiaT-catalyzed reactions

$$r_j = \frac{k_{f,j}\,[E_j]\,[\mathrm{NSD}_z]\,[\mathrm{OS}_i]}{K_{d,z}^{\mathrm{enz}}\,K_{d,i}^{\mathrm{enz}}\left(1+\frac{[\mathrm{NSD}_z]}{K_{d,z}^{\mathrm{enz}}}+\sum_{k=1}^{N_\mathrm{C}}\frac{[\mathrm{OS}_k]}{K_{d,k}^{\mathrm{enz}}}+\frac{[\mathrm{NSD}_z]}{K_{d,z}^{\mathrm{enz}}}\sum_{\substack{k=1\\k\neq i+1}}^{N_\mathrm{C}}\frac{[\mathrm{OS}_k]}{K_{d,k}^{\mathrm{enz}}}+\frac{[\mathrm{Nuc}_n]}{K_{d,n}^{\mathrm{enz}}}+\frac{[\mathrm{Nuc}_n]}{K_{d,n}^{\mathrm{enz}}}\frac{[\mathrm{OS}_{i+1}]}{K_{d,i+1}^{\mathrm{enz}}}\right)} \tag{44}$$

In these expressions, $\mathrm{OS}_i$ and $\mathrm{OS}_{i+1}$ are the reactant and product of reaction $j$, respectively; the set $\{\mathrm{OS}_k \mid k = 1, 2, \ldots, N_\mathrm{C}\}$ contains every oligosaccharide that can bind to the enzyme $E_j$; $K_{d,i}^{\mathrm{enz}}$ and $K_{d,k}^{\mathrm{enz}}$ ($\mathrm{mmol\,L^{-1}}$) are the dissociation constants for $\mathrm{OS}_i$ and $\mathrm{OS}_k$ bound to the enzyme enz $\in$ {Man I, Man II, GnT I, GnT II, GalT, FucT, SiaT}, respectively; $K_{d,z}^{\mathrm{enz}}$ and $K_{d,n}^{\mathrm{enz}}$ ($\mathrm{mmol\,L^{-1}}$) are the dissociation constants for nucleotide sugar donor $\mathrm{NSD}_z$ and nucleotide $\mathrm{Nuc}_n$ bound to the enzyme enz, where $n \in \{1, 2, 3, 4\}$; $K_{d,\mathrm{Mn}}^{\mathrm{enz}}$ refers to the dissociation constants for manganese bound to enzymes enz $\in$ {GnT I, GnT II, GalT}; $[\cdot]$ denotes the relevant intracellular concentration ($\mathrm{mmol\,L^{-1}}$); and $k_{f,j}$ $\left(\mathrm{min}^{-1}\right)$ is the kinetic rate constant of reaction $j$.

In (43)–(44), the NSD concentration inside the Golgi apparatus, $[\mathrm{NSD}_z]$, is assumed to equal 20 times of of the corresponding cytosolic concentration, $[\mathrm{NSD}_i]$, computed from the NSD synthesis model [21]. The nucleotide concentration inside the Golgi apparatus, $[\mathrm{Nuc}_n]$, is taken to be saturating and is therefore fixed at the constant values reported in Table S1 [21].

The catalytic rate constant $k_{f,j}$ follows the pH-dependent expression [45]:

$$k_{f,j} = k_{f,j}^{\max}\,\exp\left(-\frac{1}{2}\left(\frac{\mathrm{pH}^{\mathrm{golgi}}-\mathrm{pH}_{\mathrm{opt}}^{\mathrm{golgi}}}{\omega_{f,j}}\right)^2\right), \tag{45}$$

with the intra-Golgi pH estimated from the ammonia-based buffer relation

$$\mathrm{pH}^{\mathrm{golgi}} = pK_A^{\mathrm{golgi}} + \log\left(\frac{[\mathrm{Amm}]}{N_A^{\mathrm{golgi}} - [\mathrm{Amm}]}\right), \tag{46}$$

where $k_{f,j}^{\max}$ $\left(\mathrm{min}^{-1}\right)$ is the maximal kinetic rate, $\mathrm{pH}_{\mathrm{opt}}^{\mathrm{golgi}}$ is the optimal pH for enzyme $j$, $\omega_{f,j}$ is a dimensionless width parameter, $pK_A^{\mathrm{golgi}}$ is the apparent acid dissociation parameter, and $N_A^{\mathrm{golgi}}$ is a fitted buffer capacity.

Each enzyme concentration along the Golgi cisternae follows a Gaussian (bell-shaped) profile:

$$\left[E_j\right](z) = E_{j,\max}\ \exp\left[-\tfrac{1}{2}\left(\frac{z - z_{j,\max}}{\sigma_j}\right)^2\right], \tag{47}$$

where $E_{j,\max}$ (mmol L$^{-1}$) is the peak concentration, $z_{j,\max}$ is the axial position of this peak (scaled 0–1 from cis to trans), and $\sigma_j$ is the standard deviation that characterizes enzyme dispersion along the Golgi stack.

## 3. Control problem

We focus on economic control of a fed-batch bioreactor—the mode most employed in industry because of its operational flexibility. The culture runs for a total duration $T$. Nutrient supplements are added for a short period $\Delta t_{\mathrm{feed}}$ once every sampling interval $\tau = T/N$, where $N$ is the number of control intervals. No harvest stream is withdrawn until the batch ends at $t = T$.

The goal is to maximize the harvest titer of the target glycoform while satisfying all process and product constraints. A widely used performance metric is the galactosylation index (GI) at harvest,

$$\mathrm{GI}(T) = [\mathrm{FA2G1}](T) + 2\,[\mathrm{FA2G2}](T), \tag{48}$$

where the bracketed terms denote the extracellular concentrations of the indicated glycoforms. Manipulated variables $u$ (feed flow rates) and state variables $x$ (cell density, metabolite concentrations, glycoform percentages, etc.) are bounded by

$$u_{\mathrm{lb}} \le u \le u_{\mathrm{ub}}, \qquad x_{\mathrm{lb}} \le x \le x_{\mathrm{ub}}. \tag{49}$$

The manipulated inputs are the feed stream flow rates $F_s$ for $s = 1, \dots, S$. Up to three feeds can be employed—supplement medium, galactose solution, and uridine solution—so $0 \le S \le 3$. The composition of each stream is listed in Table 1.

In an ideal scenario—perfect model, exact initial state, and no disturbances—the economic optimization of the fed-batch process reduces to an offline dynamic optimization (DO-offline),

**Table 1**
Compositions of the feed streams (all concentrations in mmol L$^{-1}$).

| Metabolite | Medium | Gal solution | Urd solution |
|---|---|---|---|
| Glc | 144.37 | 0 | 0 |
| Gln | 0 | 0 | 0 |
| Lac | 0 | 0 | 0 |
| Amm | 0.06 | 0 | 0 |
| Glu | 12.19 | 0 | 0 |
| Asn | 26.99 | 0 | 0 |
| Asp | 51.95 | 0 | 0 |
| Gal | 0 | 3600 | 0 |
| Urd | 0 | 0 | 2000 |

**(DO-offline)**

$$\max_{x(t),\,u(t)} \quad \mathrm{GI}(T) \tag{50}$$

$$\text{s.t.} \quad x(t) = \mathrm{glyco}\big(u(t), p\big), \quad t \in [0, T], \tag{51}$$

$$x(0) = x_0, \tag{52}$$

$$x_{\mathrm{lb}}^{\mathrm{path}} \le x(t) \le x_{\mathrm{ub}}^{\mathrm{path}}, \quad t \in (0, T), \tag{53}$$

$$x_{\mathrm{lb}}^{\mathrm{terminal}} \le x(T) \le x_{\mathrm{ub}}^{\mathrm{terminal}}, \tag{54}$$

$$u_{\mathrm{lb}}(t) \le u(t) \le u_{\mathrm{ub}}(t), \quad t \in [0, T]. \tag{55}$$

where glyco($u, p$) denotes the implicit PDAE model (1)–(47) parameterized by $p$ (the glycosylation model parameters); $x_0$ represents the initial values of the state variables; $x^{\mathrm{path}}$ and $x^{\mathrm{terminal}}$ are path and terminal constraint bounds; and the input bounds $u_{\mathrm{lb}}(t)$ and $u_{\mathrm{ub}}(t)$ are time-dependent, allowing, for example, zero flow rates during hold periods.

In practice, the ideal assumptions rarely hold. We address model–plant mismatch by embedding parameter adaptation within the NMPC framework (ANMPC), while the remaining assumptions hold given careful pre-start measurements and carefully controlled bolus feeding.

To evaluate the closed-loop performance of ANMPC, we inject measurement noise into simulated data using the standard deviations (stds) computed from the dataset of Ref. [21]. Two measurement scenarios are considered:

(1) **Full measurement** — all state variables are available;

(2) **Industry-realistic measurement** — NSDs and certain metabolites (Asn, Asp) and glycoforms (HM, FA1G1, SIA) remain unobserved, reflecting current analytical limitations.

The potential measurement variables and their standard deviations are summarized in Tables 2–4.

## 4. Adaptive NMPC (ANMPC) framework

The adaptive NMPC loop iteratively updates model parameters, state estimates and future control actions using online measurements. Algorithm 1 summarizes the procedure:

**Table 2**
Measurement variables and standard deviations in the cell culture model.

| Variables | Std | Variables | Std |
|---|---|---|---|
| $V$ (L) | 0.001 | $c_{\mathrm{Glu}}$ (mmol L$^{-1}$) | 0.22 |
| $X$ (cells mL$^{-1}$) | 494.66 | $c_{\mathrm{Gal}}$ (mmol L$^{-1}$) | 3.80 |
| $X_{\mathrm{v}}$ (cells mL$^{-1}$) | 494.66 | $c_{\mathrm{Urd}}$ (mmol L$^{-1}$) | 0.87 |
| $c_{\mathrm{Glc}}$ (mmol L$^{-1}$) | 1.45 | $c_{\mathrm{Asn}}$ (mmol L$^{-1}$) | 0.22 |
| $c_{\mathrm{Gln}}$ (mmol L$^{-1}$) | 0.11 | $c_{\mathrm{Asp}}$ (mmol L$^{-1}$) | 0.84 |
| $c_{\mathrm{Lac}}$ (mmol L$^{-1}$) | 0.48 | $c_{\mathrm{mAb}}$ (mmol L$^{-1}$) | 11.05 |
| $c_{\mathrm{Amm}}$ (mmol L$^{-1}$) | 0.17 | | |

**Table 3**
Measurement variables and standard deviations in the NSD model (all concentrations/Stds are in mmol L$^{-1}$).

| Variables | Std | Variables | Std |
|---|---|---|---|
| $c_{\mathrm{UDP-GlcNAc}}$ | 0.25 | $c_{\mathrm{UDP-Glc}}$ | 0.22 |
| $c_{\mathrm{GDP-Fuc}}$ | 0.004 | $c_{\mathrm{GDP-Man}}$ | 0.003 |
| $c_{\mathrm{UDP-Gal}}$ | 0.18 | $c_{\mathrm{UDP-GalNAc}}$ | 0.03 |
| $c_{\mathrm{CMP-Neu5Ac}}$ | 0.005 | | |

**Table 4**
Measurement variables and standard deviations in the Golgi model (all variables/Stds are in %).

| Variables | Std | Variables | Std |
|---|---|---|---|
| $y^{\mathrm{extra}}_{\mathrm{HM}}$ | 0.64 | $y^{\mathrm{extra}}_{\mathrm{SIA}}$ | 0.12 |
| $y^{\mathrm{extra}}_{\mathrm{FA1G1}}$ | 0.15 | $y^{\mathrm{extra}}_{\mathrm{G0}}$ | 0.42 |
| $y^{\mathrm{extra}}_{\mathrm{FA2G0}}$ | 1.05 | $y^{\mathrm{extra}}_{\mathrm{G2}}$ | 0.49 |
| $y^{\mathrm{extra}}_{\mathrm{FA2G1}}$ | 0.99 | $y^{\mathrm{extra}}_{\mathrm{Man5}}$ | 0.90 |
| $y^{\mathrm{extra}}_{\mathrm{FA2G2}}$ | 0.82 | | |

**Algorithm 1 (ANMPC algorithm).**

**Initialization**
Set counter $k = 0$, current time $t = 0$. Specify: total batch duration $T$; number of control intervals $N$ and length $\tau = T/N$; feed window $\Delta t_{\mathrm{feed}}$; preparation time $\Delta t_{\mathrm{prep}}$. Times: $t_j = j\,\tau$ $(j = 0, \dots, N)$, $t_j^{\mathrm{end}} = j\,\tau + \Delta t_{\mathrm{feed}}$ $(j = 0, \dots, N-1)$, $t_j^{\mathrm{sample}} = j\,\tau - \Delta t_{\mathrm{prep}}$ $(j = 1, \dots, N)$. Initialize parameters $\hat{p}(0)$.

**Step 1: Parameter estimation**
If $k \geq 1$ and the $k$th sample has been analyzed, update $\hat{p}(k)$ using all the available measurements $\{y(j)\}_{j=1}^{k}$; otherwise skip to Step 3.

**Step 2: State estimation**
Simulate the process model from $t = 0$ to $t_k$ with $p = \hat{p}(k)$ and obtain the current state estimate $\hat{x}(k)$.

**Step 3: Dynamic optimization**
Solve the dynamic optimization problem (soft-DO-$k$) over $t_k \leq t \leq T$ with parameters $\hat{p}(k)$ and initial condition $\hat{x}(k)$, and generate the optimal control sequence $\mathbf{u} = \left[u(k)^{\top}, u(k+1)^{\top}, \dots, u(N-1)^{\top}\right]^{\top}$.

**Step 4: Bioreactor operation & sampling**
Implement $u(k)$ over $t_k \leq t \leq t_k^{\mathrm{end}}$, and acquire a new sample at $t_{k+1}^{\mathrm{sample}}$ with the analysis results of $y(k+1)$. If $t_{k+1} < T$, set $k \leftarrow k+1$ and return to Step 1; otherwise proceed to Step 5.

**Step 5: Harvest**
Shut down the bioreactor and harvest the product.

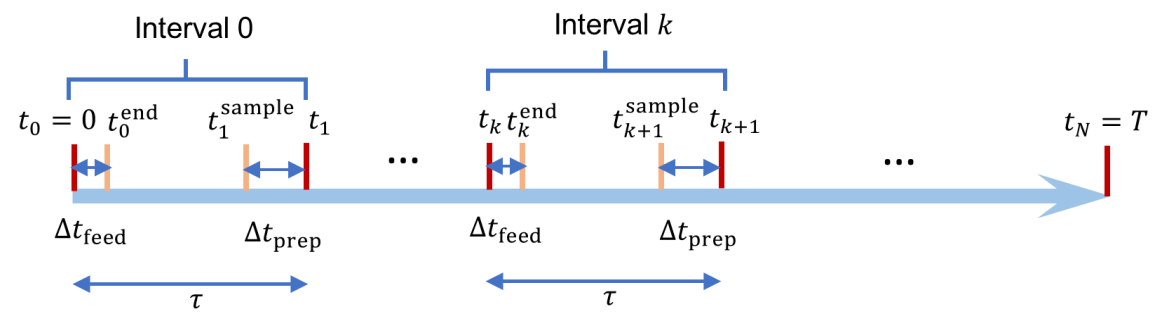

**Figure 4:** Timeline of the ANMPC algorithm.

The timeline of the ANMPC algorithm is shown in Figure 4. The preparation time $\Delta t_{\mathrm{prep}}$ in ANMPC must exceed the combined duration of analytics, parameter estimation, state estimation, and optimization to guarantee real-time feasibility. This relies on the efficient optimization algorithm introduced later. For the case studies considered here, the algorithm parameters are listed in Table 5.

**Table 5**
Parameters in the ANMPC algorithm for the glycosylation process control problem.

| Parameter | Description | Value |
|---|---|---|
| $T$ (h) | Total culture time | 288 |
| $V_0$ (L) | Initial working volume | 1.5 |
| $\Delta t_{\mathrm{feed}}$ (h) | Feed window length | 0.01 |
| $\Delta t_{\mathrm{prep}}$ (h) | Preparation time | 2 in Sections 5.2–5.3; 4 in Sections 5.4–5.5 |
| $\tau$ (h) | Length of each control interval | 24 |
| $N$ (–) | Number of control intervals | 12 |

The next subsections detail the dynamic optimization, parameter estimation, and state estimation components of the ANMPC.

## 4.1. Dynamic optimization formulation

To generate the control actions at $t_k$ in Step 3 of ANMPC, the controller solves the finite-horizon optimization (DO-$k$):

**(DO-$k$)**

$$\max_{x(t),\,u(t)} \quad \mathrm{GI}(T) \tag{50}$$

$$\text{s.t.} \quad x(t) = \mathrm{glyco}\big(u(t), \hat{p}(k)\big), \quad t \in [0, T], \tag{56}$$

$$x(k) = \hat{x}(k), \tag{57}$$

$$x_{\mathrm{lb}}^{\mathrm{path}} \leq x(t) \leq x_{\mathrm{ub}}^{\mathrm{path}}, \quad t \in (0, T), \tag{58}$$

$$x_{\mathrm{lb}}^{\mathrm{terminal}} \leq x(T) \leq x_{\mathrm{ub}}^{\mathrm{terminal}}, \tag{59}$$

$$u_{\mathrm{lb}}(t) \leq u(t) \leq u_{\mathrm{ub}}(t), \quad t \in [t_k, T]. \tag{60}$$

Relative to the offline problem (DO-offline), the model parameters $p$ are replaced by their current estimates $\hat{p}(k)$ in the problem (DO-$k$), the initial conditions are substituted by the latest state estimate $\hat{x}(k)$, and the horizon is shortened to $[t_k, T]$ because of the shrinking horizon in the fed-batch operation [6, 32]. It is expected that as data accumulate, $\hat{p}(k)$ and $\hat{x}(k)$ move closer to the ground-truth values, the optimization problem (DO-$k$) will yield progressively better predictions than the static offline optimization that uses $\hat{p}(0)$.

**Table 6**
Constraint bounds, types, and penalty weights.

| Variable | Description | Constraint type | Lower bound | Upper bound | Penalty weight |
|---|---|---|---|---|---|
| viability (%) | Cell viability | Terminal, soft | 60 | – | 10 |
| $Y^{\text{extra}}_{\text{Man5}}$ (%) | Extracellular Man5 percentage | Terminal, soft | – | 5 | 100 |
| $V$ (L) | Working volume | Terminal, hard | 0.75 | 2.25 | – |
| $F_s$ (L h$^{-1}$) $s = 1, 2, 3$ | Feed flow rate | Input | 0 | 100 | – |

Because parameter uncertainty can cause previously computed inputs to be overly aggressive, DO-$k$ may become infeasible even when DO-0 was feasible. To reduce infeasibility in the optimization, we employ a soft-constraint formulation [37],

**soft-DO-$k$:**

$$\max_{\substack{x(t),\, u(t),\\ \overline{w},\, \underline{w},\, \overline{v},\, \underline{v}}} \mathrm{GI}(T) - \int_{t_k}^{T} \rho_w^\top \big(\overline{w} + \underline{w}\big) dt - \rho_v^\top \big(\overline{v} + \underline{v}\big) \tag{61}$$

$$\begin{aligned}
\text{s.t.}\quad & \text{Eqs. (56)–(57)} \\
& x_{\text{lb}}^{\text{hard,path}} \le x^{\text{hard}}(t) \le x_{\text{ub}}^{\text{hard,path}}, \quad t \in (0, T) && (62)\\
& x_{\text{lb}}^{\text{hard,terminal}} \le x^{\text{hard}}(T) \le x_{\text{ub}}^{\text{hard,terminal}} && (63)\\
& x_{\text{lb}}^{\text{soft,path}} - \underline{w} \le x^{\text{soft}}(t) \le x_{\text{ub}}^{\text{soft,path}} + \overline{w}, \quad t \in (0, T) && (64)\\
& x_{\text{lb}}^{\text{soft,terminal}} - \underline{v} \le x^{\text{soft}}(T) \le x_{\text{ub}}^{\text{soft,terminal}} + \overline{v} && (65)\\
& \overline{w},\ \underline{w},\ \overline{v},\ \underline{v} \ge 0 && (66)\\
& u_{\text{lb}}(t) \le u(t) \le u_{\text{ub}}(t), \quad t \in [t_k, T] && (60)
\end{aligned}$$

Here the state constraints are partitioned into hard constraints, which must always be respected, and soft constraints, which may be violated at the cost of penalty weights $\rho_w$ (path) and $\rho_v$ (terminal). The variables $\overline{w}(t)$ and $\underline{w}(t)$ are slack variables for the soft path constraints, whereas $\overline{v}$ and $\underline{v}$ are slack variables for the soft terminal constraints. Problem (soft-DO-$k$) is solved only if the original problem (DO-$k$) is infeasible.

Table 6 summarizes the constrained variables along with their permissible ranges and associated penalty weights.

## 4.2. Parameter estimation formulation

The multiscale model contains ∼100 kinetic and transport parameters—far too many to calibrate accurately from the limited data available in most cases. Consequently, the parameters are adapted online using every measurement collected to date, i.e., Step 1 in ANMPC. Furthermore, the simulation with the updated parameters is used for state estimation, i.e., Step 2 in ANMPC.

We explore two windowing strategies for the online parameter estimation problem—expanding horizon estimation (EHE) and full horizon estimation (FHE). In the fed-batch process, the former is more natural as the horizon with available data becomes longer over time, so we are interested in the expanded horizon. In EHE, once the $k$th sample at $t_k^{\text{sample}}$ ($k \ge 1$) is available, we solve the problem

**EHE-$k$:**

$$\begin{aligned}
\min_{\hat{p}(k),\, x(t)} \ & \sum_{j=1}^{k} \big(y(j) - \hat{y}(j)\big)^\top V_\epsilon(j)^{-1} \big(y(j) - \hat{y}(j)\big) && (67)\\
& + \big(\hat{p}(k) - \hat{p}(0)\big)^\top P \big(\hat{p}(k) - \hat{p}(0)\big) \\
\text{s.t.}\ & x(t) = \mathrm{glyco}\big(u(t), \hat{p}(k)\big), \quad t \in \Big[0,\, t_k^{\text{sample}}\Big], && (68)\\
& x(0) = x_0, && (69)\\
& p_{\text{lb}} \le \hat{p}(k) \le p_{\text{ub}}. && (70)
\end{aligned}$$

where $y(j)$ and $\hat{y}(j)$ are the measured and predicted outputs, $V_\epsilon(j)$ is the measurement-noise covariance, and $P$ penalizes deviation from the prior parameter vector $\hat{p}(0)$. When $P = V_p^{-1}$ (the prior covariance matrix), we obtain a maximum a posteriori (MAP) estimator (e.g., [10]); when $P = 0$, we get a maximum likelihood estimator (MLE) (e.g., [3]).

EHE may overfit early data and predict physically impossible trajectories (e.g., negative concentrations) in the later horizon, jeopardizing feasibility of the subsequent NMPC problem. To reduce overfitting, we instead fit the parameters over the entire horizon $[0, T]$ even though data exist only up to $t_k^{\text{sample}}$, i.e., full-horizon estimation (FHE). In this way, parameters causing infeasible simulation beyond the measurement window are automatically discarded when solving the parameter estimation problem iteratively. To further protect against unphysical states, upper and lower bounds on the states can be incorporated into the estimation problem according to domain knowledge, resulting in the constrained FHE estimator

**cons-FHE-$k$:**

$$\begin{aligned}
\min_{\hat{p}(k),\, x(t)} \ & \text{Eq. (67)} \\
\text{s.t.}\ & x(t) = \mathrm{glyco}\big(u(t), \hat{p}(k)\big), \quad t \in [0,\, T], && (71)\\
& \text{Eqs. (69) and (70)}, \\
& x_{\text{lb}} \le x(t) \le x_{\text{ub}}, \quad t \in [0,\, T]. && (72)
\end{aligned}$$

There are two points to mention for the implementation of the online estimator:

(1) In the online parameter adaptation step, distinct estimators are applied to different submodels. Parameters in the cell culture submodel are updated with a MAP estimator $\left(P = V_p^{-1} \text{ in Eq. (67)}\right)$, whereas parameters in the NSD and Golgi submodels are fitted together by MLE $\left(P = 0 \text{ in Eq. (67)}\right)$. Alternative combinations were tested but proved inferior: applying MLE to the cell culture block led to severe overfitting in the early phase of the run, which in turn drove the simulation and subsequent optimization to diverge; conversely, employing MAP for the NSD and Golgi blocks introduced an overly strong prior that resulted in underfitting and degraded ANMPC performance. Fitting the NSD and Golgi models separately also produced underfitting and was therefore abandoned.

(2) Because the parameters span several orders of magnitude, the optimization is performed in log-space for self-normalization; that is, we estimate $\ln \hat{p}(k)$ rather than $\hat{p}(k)$. Without this approach, parameter estimation often fails due to divergence in the sensitivity-equation evaluations, where the Jacobian becomes highly ill-conditioned.

### 4.3. Solution of the optimization problems

The above dynamic optimization and parameter estimation both require repeated solutions of optimization problems subject to the multiscale glycosylation model. We adopt the control vector parameterization (CVP) method for the optimizations [22, 44], which is more reliable than simultaneous methods for strongly nonlinear dynamics–constrained problems [13]. In CVP, however, the large PDAE system described in Section 2 and its sensitivity equations must be solved tens or hundreds of times with an implicit time-stepping scheme, which dominates the computational cost.

To make the optimization tractable, we accelerate each model evaluation with the parallel QSS algorithm introduced in our earlier work, which reduces overall simulation/optimization time by two orders of magnitude [26]. The idea relies on a clear time-scale separation: the Golgi reactions equilibrate much faster than the outer-layer dynamics. Long before extracellular conditions change appreciably, the Golgi has already reached a local steady state. Under the steady state, the Golgi model becomes a DAE (along the length of the Golgi apparatus) instead of the original PDAE, simplifying the computation significantly. The complete QSS algorithm is described briefly as Algorithm 2.

**Algorithm 2 (QSS algorithm).**

**Initialization**. Choose a set of time points $\mathcal{T} = \{t_s\}_{s=1}^{N_{\text{QSS}}}$.

**Step 1: Get env**. Integrate the cell-culture model (omitting glycoprotein balances) and the NSD model to obtain the environmental states $\text{env}(t_s)$ at each $t_s$.
**Step 2: Get $Y_i^{\textbf{intra}}$**. For every $t_s$, solve the steady-state Golgi DAE with $\text{env}(t_s)$ as inputs, yielding intracellular glycoform fractions $Y_i^{\text{intra}}(t_s)$. The $N_{\text{QSS}}$ Golgi simulations are fully multithread parallelizable.
**Step 3: Get $Y_i^{\textbf{extra}}$**. Reintegrate the complete cell culture model to produce the extracellular glycoform trajectories $Y_i^{\text{extra}}(t_s)$, where all the $Y_i^{\text{intra}}(t_s)$ in the model are treated as time-varying parameters and are obtained from Step 2.

Following Ref. [26], we achieve satisfactory accuracy for $Y_i^{\text{extra}}$ when the QSS time grid $\mathcal{T}$ includes all event instants (feed start/stop, sampling, etc.) and is augmented with 100 uniformly spaced points over $[0, T]$.

## 5. Case study

To evaluate the proposed ANMPC framework, we first calibrate the multiscale model to the dataset of Ref. [21], obtaining a MAP parameter vector $\ln p^*$ and its covariance $V_{\ln(p)}$. We treat $\ln p^*$ as the ground truth. Five distinct initial guesses $\ln \hat{p}(0)$ are then generated on the boundary of the 50% confidence ellipsoid defined by $V_{\ln(p)}$ so that the conclusions based on the case studies are more reliable. These $\ln \hat{p}(0)$ vectors, together with $V_{\ln(p)}$, serve as the prior mean and covariance in the online estimation step. The ground-truth and initial parameter values are listed in Tables S2–S4, while the covariances for submodels are provided in the Excel file in Supplementary 2. For convenience, we denote the five initial parameter sets as `Init_0`, `Init_1`, `Init_2`, `Init_3`, and `Init_4`.

All computations are performed on a Windows 11 machine with a 12th-Generation Intel Core i7-12700H CPU (14 cores/20 threads, 2.30 GHz) and 16 GB RAM. Python 3.9 [43] and CasADi 3.6.5 [1] were used as the computation platform. The nonlinear programs (NLPs) arising from the CVP formulation are solved with PySQP, an in-house sequential quadratic programming (SQP) solver based on the I-SQP algorithm [25] with a watchdog-technique enhancement. The DAE systems in the QSS algorithm are integrated with the IDAS solver [8, 12] via CasADi, using 15 threads for parallel computation. Under these settings, the total runtime per control move (parameter estimation + state estimation + dynamic optimization) is less than 2 h in all case studies. In Sections 5.2 and 5.3, we assume rapid analytics and set the preparation time $\Delta t_{\text{prep}} = 2$ h; in Sections 5.4 and 5.5, we use a more realistic $\Delta t_{\text{prep}} = 4$ h to account for assay time.

### 5.1. Open-loop optimization study

Because supplement medium is routinely fed whereas Gal and Urd are not, we first solve the open-loop optimization with different combinations of feed streams to study the influence of Gal and Urd feeding. Here, the ground-truth parameters are used, and the optimization results are shown in Table 7.

**Table 7**
Optimized fed-batch performance with different feeds.

| Feed streams | Harvest volume (L) | GI (mg $L^{-1}$) |
|---|---|---|
| None (batch) | 1.38 | 174.27 |
| Medium only | 2.18 | 247.28 |
| Medium + Gal | 2.23 | 404.00 |
| Medium + Urd | 2.23 | 405.25 |
| Medium + Gal + Urd | 2.23 | 405.25 |

Adding only the supplement medium raises the galactosylation index (GI) by about 42% relative to the batch baseline. Supplementing the medium with either galactose or uridine produces a much larger benefit: the optimal GI is roughly 64% higher than with medium alone (and more than double the batch value). Adding uridine gives a slightly (0.3%) higher GI than adding galactose. Introducing both galactose and uridine simultaneously provides no further gain, because the optimizer drives the galactose feed rate to zero whenever uridine is available. Thus, only one of the two supplemental streams is required; the decision between galactose and uridine should be based on robustness and operational practicality rather than peak GI alone.

Figure 5 compares the optimal time courses obtained with different feed strategies. Figure 5a confirms the ranking in Table 7: every feeding strategy raises the harvest GI relative to the batch baseline, but the trajectories diverge markedly once galactose or uridine is introduced. When only supplement medium is added, the improvement in extracellular GI stems almost entirely from biomass growth—the viable cell density roughly doubles (Figure 5b)—even though the intracellular GI falls slightly (Figure 5d). Adding galactose or uridine in addition to the medium changes the picture: both supplements elevate the intracellular GI by roughly 20% (Figure 5d), on top of the higher viable cell density achieved with medium feeding, yielding the pronounced increase in extracellular GI shown in Figure 5a. The increase in intracellular GI is driven by an elevated UDPGal concentration, a key donor for galactosylation (Figure 5c). Once the UDPGal concentration is sufficiently high, however, the incremental benefit plateaus; this is evident from the close overlap of the green dash–dot and red dotted curves in Figure 5d, which correspond to Gal and Urd supplement scenarios, respectively.

Although supplementing the reactor with Urd yields a slightly higher GI than Gal, the process is much more sensitive to Urd dosing. We quantify the response to the initial Urd pulse volume $V_{\text{Urd}}$ by varying its value from 0 mL to 20 mL while holding the medium feed fixed. Figure 6 shows a sharp optimum at 0.014 mL, giving a harvest GI of 405 mg $L^{-1}$. Increasing the pulse to just 1 mL lowers GI to about 350 mg $L^{-1}$—a 13.6 % drop—demonstrating that small deviations from the optimum can markedly degrade product quality and, therefore, robustness in the presence of

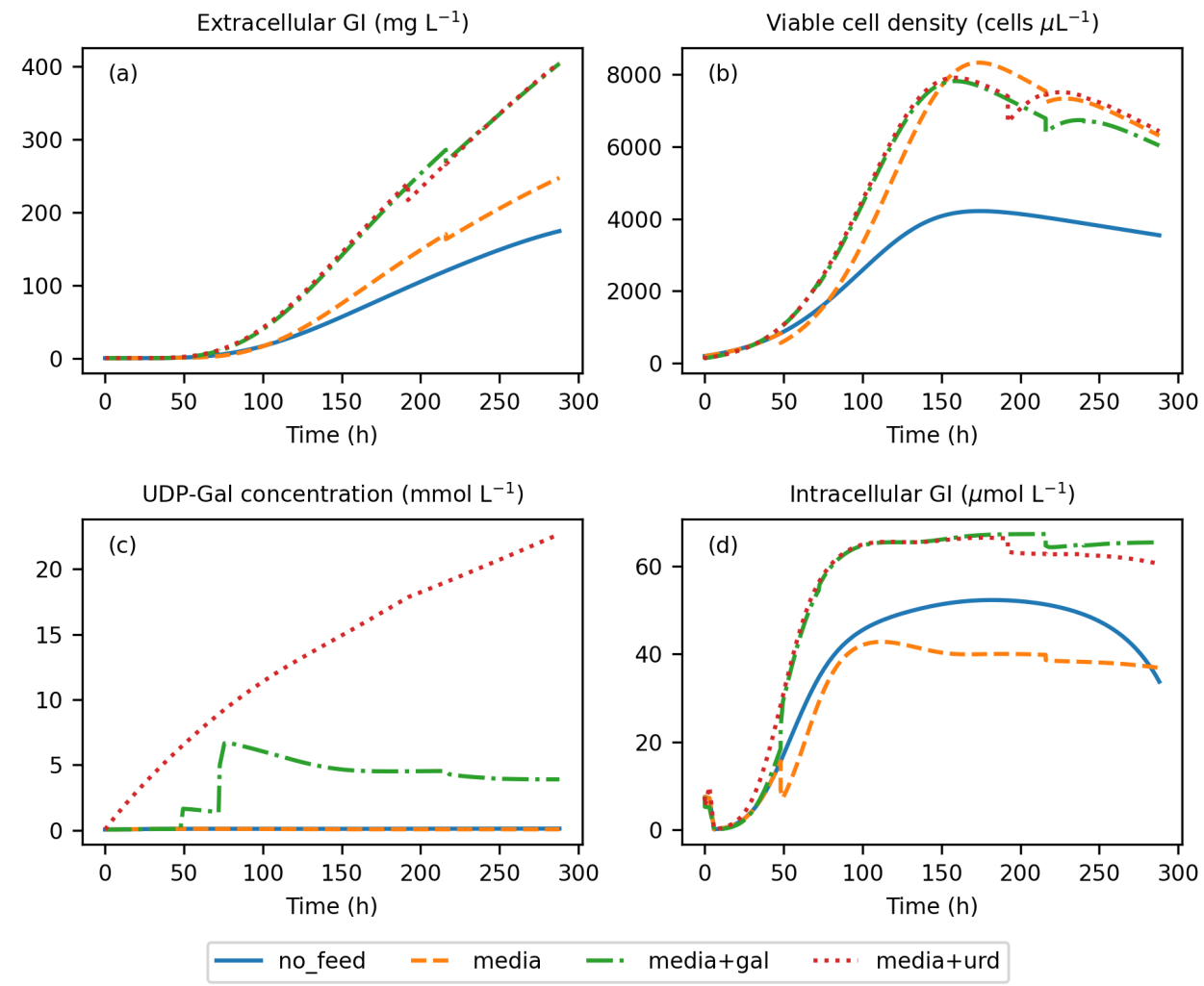

**Figure 5:** Optimal trajectories of crucial variables when using different feed streams.

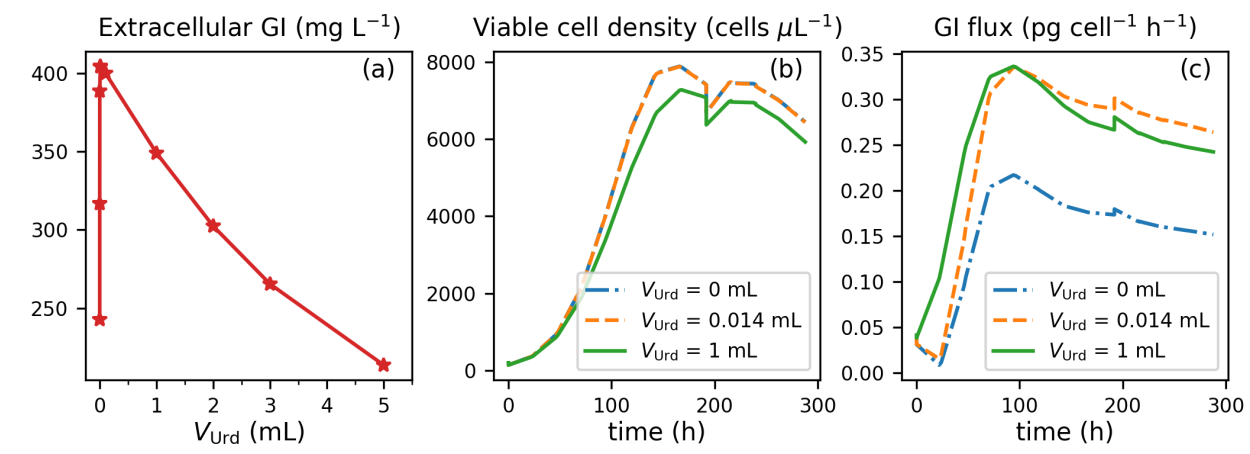

**Figure 6:** a) Sensitivity of harvest GI to the Urd feed volume on day 1, b) viable cell density trajectories for three pulse feed volumes, and c) intracellular GI trajectories for the same three different three pulse feed volumes.

model–plant mismatch or input disturbances. Figures 6bc further indicate that the GI loss is driven primarily by reduced viable cell density throughout the run and by lower intracellular GI during the high cell density phase.

In contrast, the harvest GI is far less sensitive to the galactose pulse volume $V_{\text{Gal}}$ near its optimum—15 mL on day 3 and 58 mL on day 4 , with no Gal on the other days. Figures 7a–c illustrate this robustness. When $V_{\text{Gal}}$ is above 10 mL on each of the first four days (Figure 7a), the harvest GI changes only slowly, and declines rapidly only when $V_{\text{Gal}}$ falls below 5 mL. Moreover, varying the Gal dose on a single day while keeping the other days at their optimal values (Figs. 7bc) produces only minor changes in GI. These results confirm that galactose feeding is a potent yet forgiving lever for glycosylation control, whereas uridine requires near-perfect tuning to realize its modest advantage.

Given the above analysis, the subsequent closed-loop studies employ medium + Gal as the manipulated streams and omit Urd.

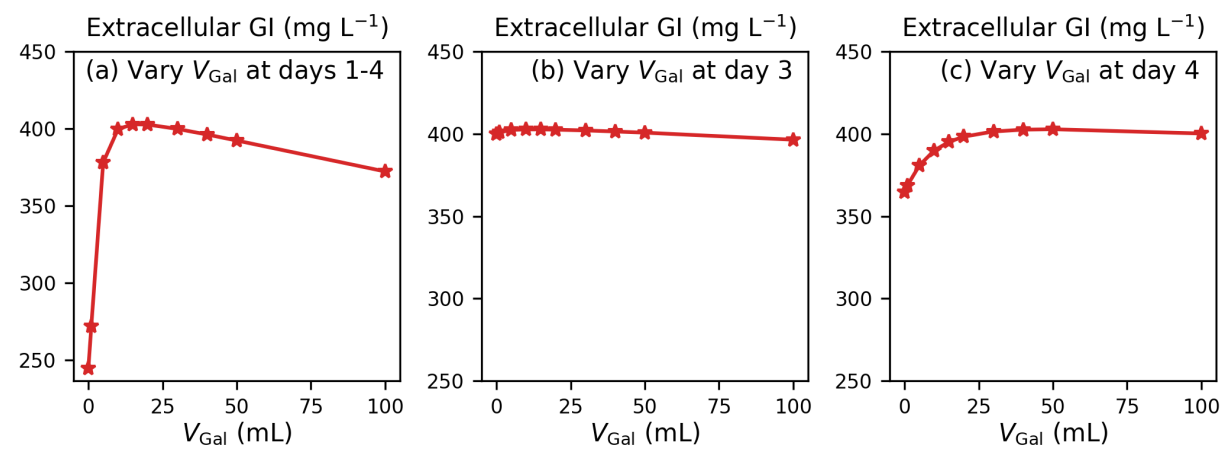

**Figure 7:** Sensitivity of harvest GI to the Gal feed volume, a) vary cap V sub Gal simultaneously on days 1–4, b) vary cap V sub Gal on day 3, c) vary cap V sub Gal on day 4.

**Table 8**
Merits under different control algorithms.

| | Init_0 | Init_1 | Init_2 | Init_3 | Init_4 |
|---|---|---|---|---|---|
| Open loop | 185.88 | 272.64 | 231.17 | 219.86 | 274.78 |
| State NMPC | 276.08 | —[a] | —[a] | 202.26 | 367.22 |
| ANMPC | 216.82 | 362.89 | 220.38 | 218.30 | 380.92 |
| Cell culture | | | 8.70 | | |
| Ground truth | | | 404.00 | | |

[a] State NMPC terminated prematurely for Init_1 and Init_2.

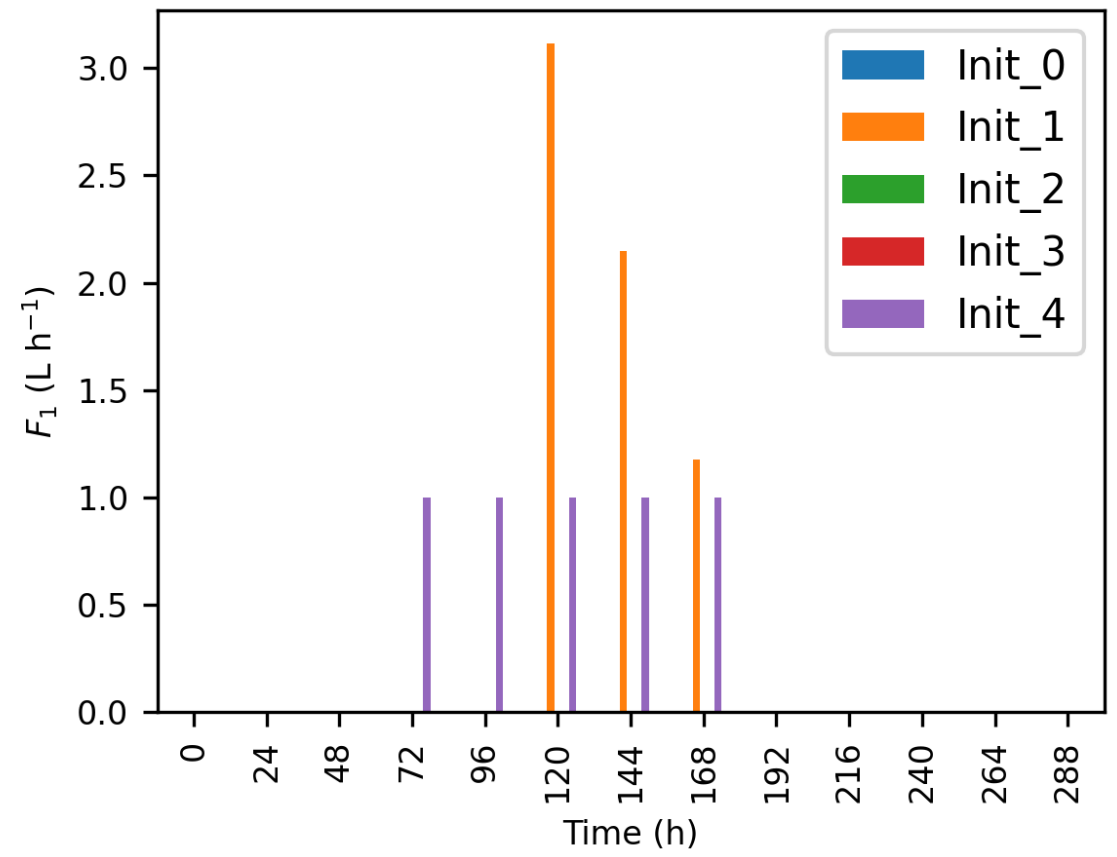

**Figure 8:** Gal feed rate profiles under ANMPC for different initial parameter sets. The feed window duration is 0.01 h.

## 5.2. Open-loop optimization vs. State NMPC vs. ANMPC

To benchmark the ANMPC scheme, we compared five controllers, all initiated from the same inaccurate parameter vector:

(1) Ground truth — dynamic optimization with the ground-truth parameters $p^*$.

(2) Cell culture — dynamic optimization with the objective of maximizing titer, which is a common scenario when there is only the cell culture model. Here, Ground-truth parameters $p^*$ are used.

(3) Open loop — dynamic optimization with the inaccurate parameters $\hat{p}(0)$.

(4) State NMPC — including feedback on measured states and using fixed parameters $\hat{p}(0)$.

(5) ANMPC — adjust parameters online with the measured states.

To compare with state NMPC more easily, we assume full state measurement in this section, while the more realistic measurement scheme will be considered in Sections 5.4 and 5.5.

Because model–plant mismatch may drive viability at harvest below its 60% limit and thereby activate the soft-constraints optimization, we evaluate each controller with a penalized merit function

$$\begin{aligned}\text{Merit} = \text{GI}(T) - \rho_{v,1}\ \max\big(0,\ 60 - \text{viability}\big)\\ - \rho_{v,2}\ \max\big(0,\ Y_{\text{Man5}}^{\text{extra}} - 5\big),\end{aligned} \tag{73}$$

using the same penalty weights $\rho_v$ as in the soft-constrained formulation (soft-DO-$k$; see Table 5). Table 8 reports the merits across the five parameter scenarios.

From Table 8, optimizing the cell culture model only will cause extremely low GI (50 times smaller than the ground-truth solution) because the operating condition of maximizing titer leads to very little FA2G2 and FA2G1 synthesized in the Golgi apparatus (less than 1% for most of the culture time). Open-loop optimization performs second worst for most initializations, producing merits roughly one-half of the ground-truth optima. ANMPC shows a substantial gain (over 30%) for Init_1 and Init_4, but offers only marginal or no improvement over open-loop optimization for the other runs. State NMPC terminates prematurely for Init_1 and Init_2 due to infeasible simulations during optimization. Overall, none of the controllers consistently achieve satisfactory performance.

A closer look at the ANMPC operation strategies shows that the low-merit cases (Init_0, Init_2, Init_3) never feed Gal, whereas the high-merit cases (Init_1, Init_4) apply appropriate Gal dosing (Figure 8). Consistent with Table 7, omitting Gal caps the ground-truth optimal GI at just 247 mg L$^{-1}$. This behavior occurs because the inaccurate initial parameters undervalue Gal's positive impact on galactosylation, so the optimizer chooses not to add Gal. Without Gal feeding, the Gal-related parameters are not excited and therefore cannot be updated, which masks the potential advantage of ANMPC.

## 5.3. Enforcing a minimum galactose dose

The previous section showed that, when ANMPC selects zero galactose feeding, the Gal-related parameters cannot be updated and performance stalls. To force at least minimal excitation, we impose a lower bound of 1 mL on the Gal feed volume during day 1—a negligible volume compared with over 1 L working volume. Table 9 compares the resulting merits for the various control strategies under this specification.

**Table 9**
Merits after imposing a 1 mL Gal feeding lower bound on day 1.

| | Init_0 | Init_1 | Init_2 | Init_3 | Init_4 |
|---|---|---|---|---|---|
| Open loop | 183.36 | 324.88 | 197.27 | 166.23 | 201.67 |
| State NMPC | 224.55 | —[a] | —[a] | 197.11 | 378.11 |
| ANMPC | 380.37 | 398.29 | 389.53 | 386.79 | 387.06 |
| Ground truth | | | 404.00 | | |

[a] State NMPC terminated prematurely for Init_1 and Init_2.

Comparing Table 8 (no Gal feed lower bound) with Table 9 (1 mL lower bound on day 1), the open-loop and state-NMPC controllers show almost no progress: their merits remain roughly half of the 404 benchmark, and state NMPC still fails to converge for two of the five initializations. By contrast, ANMPC improves markedly once the small day-1 galactose dose is enforced. Its merits lie within 6% of the ground-truth optimum for every initialization and exhibit little scatter, indicating substantially greater robustness. Relative to ANMPC without a Gal lower bound, the merit increases by more than 75% for Init_0, Init_2, and Init_3. Compared with open-loop optimization and state NMPC under the same Gal feed constraint, ANMPC delivers gains of up to 130% and 96% (Init_3), respectively. These results confirm that introducing even a minimal Gal feed early in the cell culture provides sufficient excitation to update Gal-related parameters and enables the adaptive controller to approach its full potential.

Figure 9 illustrates the input profiles and trajectories for Init_0 alongside the ground-truth policy. The ground-truth optimization splits the medium addition into two pulses—one at the start of culture and a second after nine days (Figure 9a). This staged strategy moderates the late-phase viability drop and satisfies the terminal viability limit (Figure 9c). In contrast, the other controllers deliver nearly the entire medium volume at the outset, pushing the broth quickly to the reactor's upper volume limit and leaving little capacity for a compensatory feed later; the resulting nutrient depletion leads the NMPC variants to violate the viability constraint. The second key difference appears in Figure 9b: substantial Gal supplementation occurs only in the ground-truth and ANMPC policies, yielding a markedly higher GI in Figure 9d. Like open-loop and state NMPC, ANMPC omits Gal during the first three days due to initial model mismatch. Once the parameter estimates improve, however, ANMPC initiates Gal feeding; from that point forward, its GI trajectory closely tracks the ground-truth benchmark and remains well above the curves from the other two controllers.

The contrasting control performances can be traced to how each algorithm handles model–plant mismatch, as illustrated by the evolution of the GI prediction errors at the harvest time in Figure 10. At the start of the run (day 0) the three schemes—open-loop optimization, state NMPC, and ANMPC—share the same parameter set, so their GI prediction errors are identical, about 200 mg L$^{-1}$. Thereafter, the open-loop error remains unchanged because no feedback is used to refine the model. In both state NMPC and ANMPC, the error decreases as new measurements arrive, but ANMPC's errors are consistently smaller than those of state NMPC. Within two days, the GI error under ANMPC drops below 50 mg L$^{-1}$ (except on days 3–5), whereas state NMPC does not reach this level until about day nine. The rapid improvement in ANMPC stems from its online parameter adaptation, whereas the slower drift in state NMPC reflects only the shrinking prediction horizon and the availability of full-state measurements, not improved model accuracy. By the time state NMPC finally achieves a similar error magnitude, only three days remain to influence the harvest GI, leaving ANMPC with the practical advantage.

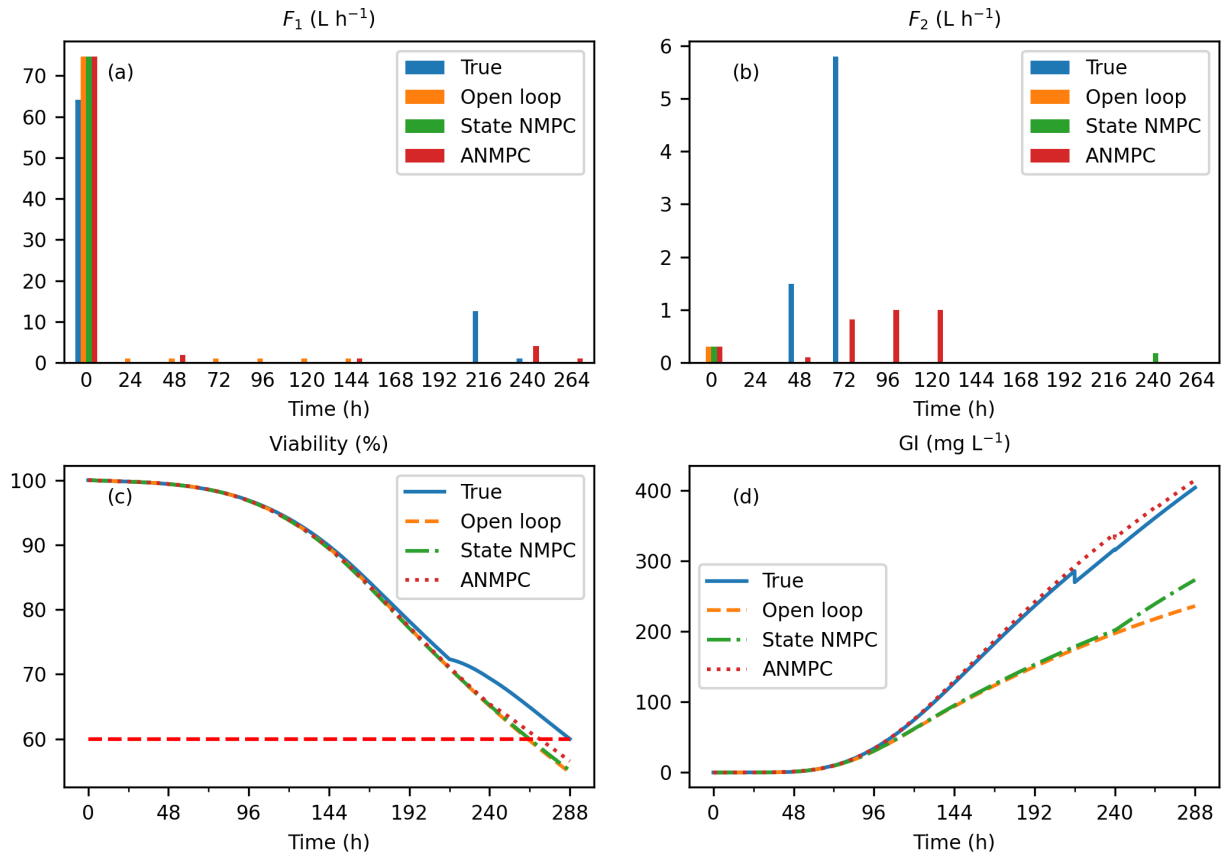

**Figure 9:** Control inputs and key state trajectories for the control schemes starting from Init_0 and the ground-truth optimization: (a) medium feed flow rate, (b) Gal feed flow rate, (c) cell viability, and (d) GI.

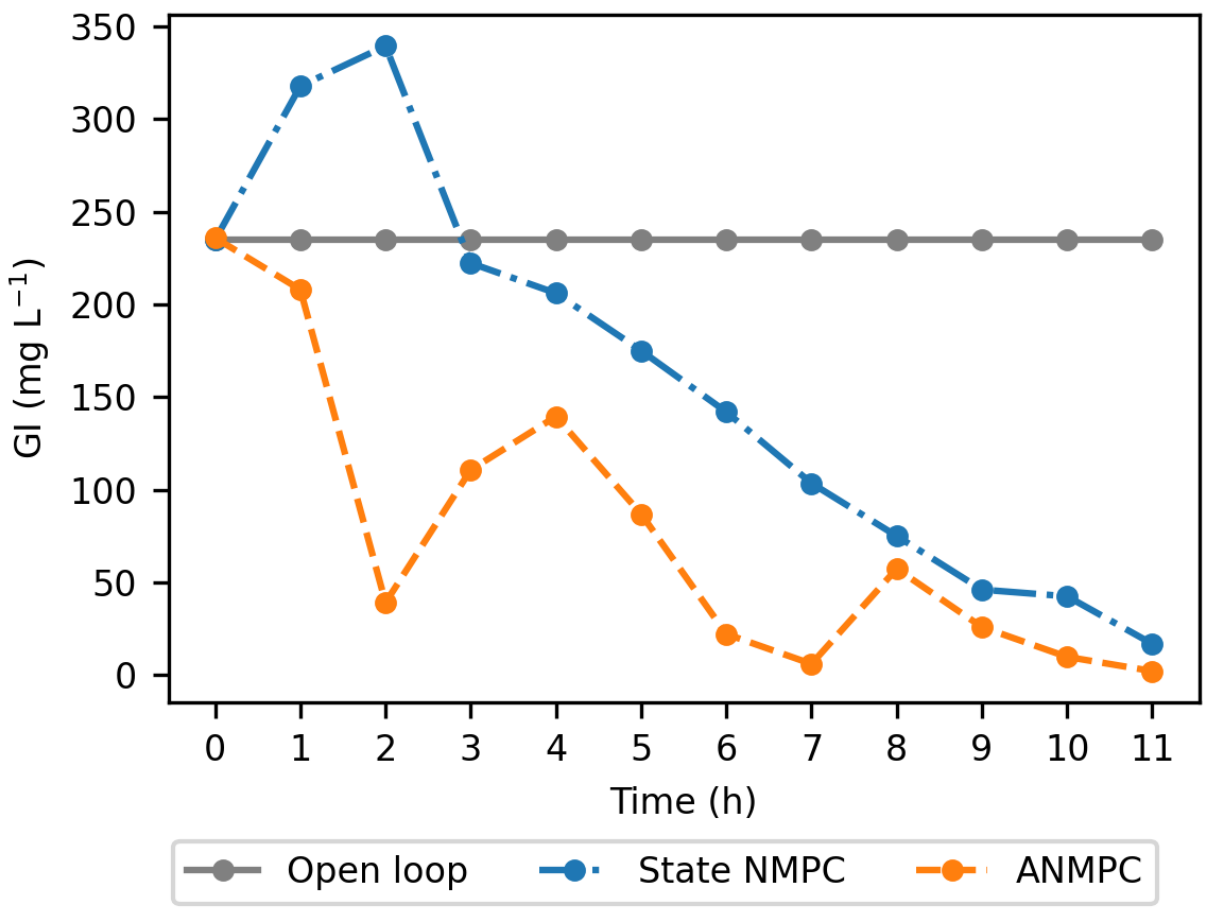

**Figure 10:** Evolution of the GI prediction error at harvest time when starting from Init_0.

## 5.4. ANMPC with fewer measurements and longer preparation time

Here we evaluate ANMPC under more realistic conditions: (1) partial measurement availability—no readings for extracellular Asn and Asp, the extracellular glycan fractions HM, FA1G1, and SIA, or any intracellular NSD concentrations; (2) a longer preparation time of 4 h to account for computation and assay turnaround. The results are summarized in Table 10.

**Table 10**
ANMPC under different measurements and preparation time.

| | Init_0 | Init_1 | Init_2 | Init_3 | Init_4 | Avg |
|---|---|---|---|---|---|---|
| Full[a] + $\Delta t_{\text{prep}} = 2$ h | 380.37 | 398.29 | 389.53 | 386.79 | 387.06 | 388.41 |
| Partial[b] + $\Delta t_{\text{prep}} = 2$ h | 387.81 | 395.60 | 392.24 | 388.52 | 384.52 | 389.74 |
| Partial[b] + $\Delta t_{\text{prep}} = 4$ h | 386.18 | 395.25 | 387.89 | 386.62 | 379.40 | 387.07 |

[a] "Full" = all variables measured.
[b] "Partial" = excludes extracellular Asn and Asp, extracellular glycan fractions (HM, FA1G1, SIA), and all intracellular NSD concentrations.

Table 10 shows that, for a given initial parameter set, the penalized merit varies by less than 2% across different measurement schemes and preparation times. In particular, using partial measurements and a 4 h preparation window does not necessarily degrade performance relative to the full measurements and 2 h baseline. Therefore, a leaner measurement set and a relaxed preparation window can be adopted as a more economical choice without sacrificing control quality.

## 6. Conclusions

We applied an ANMPC framework to regulate the N-glycosylation in a fed-batch mAb bioreactor. The controller uses a high-fidelity multiscale PDAE model that links extracellular operating conditions to intracellular Golgi reactions to predict both mAb productivity and glycan profiles. To mitigate model–plant mismatch, parameters are updated online as new measurements arrive, after which a shrinking-horizon optimization recomputes the inputs. When the nominal dynamic optimization becomes infeasible, a soft-constraint formulation preserves solvability. To prevent unphysical predictions (e.g., negative concentrations), parameter estimation is performed over the full fed-batch horizon with state bounds. Computationally, parameter estimation and dynamic optimization yield large-scale, model-constrained problems that are challenging to solve. We therefore employ the control-vector parameterization (CVP) method for reliable convergence and accelerate embedded simulations via a parallel QSS algorithm.

Sensitivity studies show that both Gal and Urd can raise GI by about 64% relative to feeding medium alone to the bioreactor, but Urd is overly dosage-sensitive for control. Across case studies with multiple initial parameter sets, ANMPC delivers up to 130% and 96% higher penalized merit than open-loop optimization and state NMPC, respectively—provided a small day-1 Gal lower bound is imposed to excite Gal-related dynamics and enable adaptation. Moreover, using fewer measurements and a longer preparation time (4 h) did not materially degrade ANMPC performance, supporting practical deployability.

Future work will explore stochastic ANMPC (e.g., chance-constrained MPC) to balance constraint satisfaction and performance under uncertainty, and experimental validation of the proposed strategies.

## Acknowledgments

This work was supported by a Project Award Agreement from the National Institute for Innovation in Manufacturing Biopharmaceuticals (NIIMBL), with financial assistance from the U.S. Department of Commerce, National Institute of Standards and Technology (NIST), awards 70NANB17H002 and 70NANB20H037. We thank Professor Cleo Kontoravdi for generously providing the experimental data reported in her publication.

## Code availability

Code will be made publicly available upon acceptance. A reviewer link can be provided to editors on request.

## References

[1] Andersson, J.A.E., Gillis, J., Horn, G., Rawlings, J.B., Diehl, M., 2019. CasADi: A software framework for nonlinear optimization and optimal control. Mathematical Programming Computation 11, 1–36. doi:10.1007/s12532-018-0139-4.

[2] Batra, J., Rathore, A.S., 2016. Glycosylation of monoclonal antibody products: Current status and future prospects. Biotechnology Progress 32, 1091–1102. doi:10.1002/btpr.2366.

[3] Beck, J.V., Arnold, K.J., 1977. Parameter Estimation in Engineering and Science. Wiley.

[4] Dan, A., Liu, B., Patil, U., Manuraj, B.N.M., Gandhi, R., Buchel, J., Chundawat, S.P.S., Guo, W., Ramachandran, R., 2025. Machine learning model-based design and model predictive control of a bioreactor for the improved production of mammalian cell-based biotherapeutics. Control Engineering Practice 156, 106198. doi:10.1016/j.conengprac.2024.106198.

[5] Dewasme, L., Fernandes, S., Amribt, Z., Santos, L.O., Bogaerts, P., Vande Wouwer, A., 2015. State estimation and predictive control of fed-batch cultures of hybridoma cells. Journal of Process Control 30, 50–57. doi:10.1016/j.jprocont.2014.12.006.

[6] Eaton, J.W., Rawlings, J.B., 1990. Feedback control of chemical processes using on-line optimization techniques. Computers & Chemical Engineering 14, 469–479. doi:10.1016/0098-1354(90)87021-G.

[7] Federici, M., Lubiniecki, A., Manikwar, P., Volkin, D.B., 2013. Analytical lessons learned from selected therapeutic protein drug comparability studies. Biologicals 41, 131–147. doi:10.1016/j.biologicals.2012.10.001.

[8] Gardner, D.J., Reynolds, D.R., Woodward, C.S., Balos, C.J., 2022. Enabling new flexibility in the SUNDIALS suite of nonlinear and

differential/algebraic equation solvers. ACM Transactions on Mathematical Software 48, 31. doi:10.1145/3539801.

[9] Goetze, A.M., Liu, Y.D., Zhang, Z., Shah, B., Lee, E., Bondarenko, P.V., Flynn, G.C., 2011. High-mannose glycans on the Fc region of therapeutic IgG antibodies increase serum clearance in humans. Glycobiology 21, 949–959. doi:10.1093/glycob/cwr027.

[10] Gunawan, R., Jung, M.Y.L., Seebauer, E.G., Braatz, R.D., 2003. Maximum a posteriori estimation of transient enhanced diffusion energetics. AIChE Journal 49, 2114–2123. doi:10.1002/aic.690490819.

[11] Hajizadeh, I., Rashid, M., Sevil, M., Brandt, R., Samadi, S., Hobbs, N., Cinar, A., 2018. Adaptive model predictive control for nonlinearity in biomedical applications, in: 6th IFAC Conference on Nonlinear Model Predictive Control (NMPC 2018), pp. 368–373. doi:10.1016/j.ifacol.2018.11.061.

[12] Hindmarsh, A.C., Brown, P.N., Grant, K.E., Lee, S.L., Serban, R., Shumaker, D.E., Woodward, C.S., 2005. SUNDIALS: Suite of nonlinear and differential/algebraic equation solvers. ACM Transactions on Mathematical Software 31, 363–396. doi:10.1145/1089014.1089020.

[13] Hong, W., Wang, S., Li, P., Wozny, G., Biegler, L.T., 2006. A quasi-sequential approach to large-scale dynamic optimization problems. AIChE Journal 52, 255–268. doi:10.1002/aic.10625.

[14] Jabarivelisdeh, B., Carius, L., Findeisen, R., Waldherr, S., 2020. Adaptive predictive control of bioprocesses with constraint-based modeling and estimation. Computers & Chemical Engineering 135, 106744. doi:10.1016/j.compchemeng.2020.106744.

[15] James, D.C., 2005. Control of recombinant monoclonal antibody effector functions by Fc N-glycan remodeling in vitro. Biotechnology Progress 21, 1644–1652. doi:10.1002/btpr.50228.

[16] Jefferis, R., 2009. Glycosylation as a strategy to improve antibody-based therapeutics. Nature Reviews Drug Discovery 8, 226–234. doi:10.1038/nrd2804.

[17] Kappatou, C.D., Ehsani, A., Niedenführ, S., Mhamdi, A., Schuppert, A., Mitsos, A., 2020. Quality-targeting dynamic optimization of monoclonal antibody production. Computers & Chemical Engineering 142, 107004. doi:10.1016/j.compchemeng.2020.107004.

[18] Kontoravdi, C., Pistikopoulos, E.N., Mantalaris, A., 2010. Systematic development of predictive mathematical models for animal cell cultures. Computers & Chemical Engineering 34, 1192–1198. doi:10.1016/j.compchemeng.2010.03.012.

[19] Kornfeld, R., Kornfeld, S., 1985. Assembly of asparagine-linked oligosaccharides. Annual Review of Biochemistry 54, 631–664. doi:10.1146/annurev.bi.54.070185.003215.

[20] Kotidis, P., 2021. Enhanced Understanding of Protein Glycosylation in CHO Cells Through Computational Tools and Experimentation. Ph.D. thesis. Imperial College London. doi:10.25560/100281.

[21] Kotidis, P., Jedrzejewski, P.M., Sou, S.N., Sellick, C., Polizzi, K., Jimenez del Val, I., Kontoravdi, C., 2019. Model-based optimization of antibody galactosylation in cho cell culture. Biotechnology and Bioengineering 116, 1612–1626. doi:10.1002/bit.26960.

[22] Kraft, D., 1985. On converting optimal control problems into nonlinear programming problems, in: Schittkowski, K. (Ed.), Computational Mathematical Programming. Springer Berlin Heidelberg, pp. 261–280.

[23] Lu, F., Toh, P.C., Burnett, I., Li, F., Hudson, T., Amanullah, A., Li, J., 2013. Automated dynamic fed-batch process and media optimization for high productivity cell culture process development. Biotechnology and Bioengineering 110, 191–205. doi:10.1002/bit.24602.

[24] Luo, Y., Kurian, V., Song, L., Wells, E.A., Robinson, A.S., Ogunnaike, B.A., 2023. Model-based control of titer and glycosylation in fed-batch mab production: Modeling and control system development. AIChE Journal 69, e18075. doi:10.1002/aic.18075.

[25] Ma, Y., Gao, X., Liu, C., Li, J., 2024. Improved sqp and slsqp algorithms for feasible path-based process optimisation. Computers & Chemical Engineering 188, 108751. doi:10.1016/j.compchemeng.2024.108751.

[26] Ma, Y., Guo, J., Maloney, A.J., Braatz, R.D., 2025. Quasi-steady-state approach for efficient multiscale simulation and optimization of mAb glycosylation in CHO cell culture. Chemical Engineering Science 318, 122162. doi:10.1016/j.ces.2025.122162.

[27] Majewska, N.I., Tejada, M.L., Betenbaugh, M.J., Agarwal, N., 2020. N-glycosylation of IgG and IgG-like recombinant therapeutic proteins: Why is it important and how can we control it? Annual Review of Chemical and Biomolecular Engineering 11, 311–338. doi:10.1146/annurev-chembioeng-102419-010001.

[28] Maloney, A.J., 2021. Case Studies in the Modeling and Control of Continuous Pharmaceutical Manufacturing Processes. Ph.D. thesis. Massachusetts Institute of Technology. URL: https://dspace.mit.edu/handle/1721.1/158315.

[29] Matsumiya, S., Yamaguchi, Y., Saito, J., Nagano, M., Sasakawa, H., Otaki, S., Satoh, M., Shitara, K., Kato, K., 2007. Structural comparison of fucosylated and nonfucosylated Fc fragments of human immunoglobulin G1. Journal of Molecular Biology 368, 767–779. doi:10.1016/j.jmb.2007.02.034.

[30] McCamish, M., Woollett, G., 2013. The continuum of comparability extends to biosimilarity: How much is enough and what clinical data are necessary? Clinical Pharmacology & Therapeutics 93, 315–317. doi:10.1038/clpt.2013.17.

[31] Mullard, A., 2021. FDA approves 100th monoclonal antibody product. Nature Reviews Drug Discovery 20, 491–495. doi:10.1038/d41573-021-00079-7.

[32] Nagy, Z.K., Braatz, R.D., 2003. Robust nonlinear model predictive control of batch processes. AIChE Journal 49, 1776–1786. doi:10.1002/aic.690490715.

[33] Pickhardt, R., 2000. Nonlinear modelling and adaptive predictive control of a solar power plant. Control Engineering Practice 8, 937–947. doi:10.1016/S0967-0661(00)00009-5.

[34] Rawlings, J.B., Mayne, D.Q., Diehl, M., 2017. Model Predictive Control: Theory, Computation, and Design. volume 2. 2 ed., Nob Hill Publishing.

[35] Sarna, S., Patel, N., Corbett, B., McCready, C., Mhaskar, P., 2023. Process-aware data-driven modelling and model predictive control of bioreactor for the production of monoclonal antibodies. The Canadian Journal of Chemical Engineering 101, 2677–2692. doi:10.1002/cjce.24752.

[36] Sauer, P.W., Burky, J.E., Wesson, M.C., Sternard, H.D., Qu, L., 2000. A high-yielding, generic fed-batch cell culture process for production of recombinant antibodies. Biotechnology and Bioengineering 67, 585–597. doi:10.1002/(SICI)1097-0290(20000305)67:5<585::AID-BIT9>3.0.CO;2-H.

[37] Scokaert, P.O.M., Rawlings, J.B., 1999. Feasibility issues in linear model predictive control. AIChE Journal 45, 1649–1659. doi:10.1002/aic.690450805.

[38] Sha, S., Agarabi, C., Brorson, K., Lee, D.Y., Yoon, S., 2016. N-glycosylation design and control of therapeutic monoclonal antibodies. Trends in Biotechnology 34, 835–846. doi:10.1016/j.tibtech.2016.02.013.

[39] Spahn, P.N., Hansen, A.H., Hansen, H.G., Arnsdorf, J., Kildegaard, H.F., Lewis, N.E., 2016. A Markov chain model for N-linked protein glycosylation – Towards a low-parameter tool for model-driven glycoengineering. Metabolic Engineering 33, 52–66. doi:10.1016/j.ymben.2015.10.007.

[40] St. Amand, M.M., Radhakrishnan, D., Robinson, A.S., Ogunnaike, B.A., 2014a. Identification of manipulated variables for a glycosylation control strategy. Biotechnology and Bioengineering 111, 1957–1970. doi:10.1002/bit.25251.

[41] St. Amand, M.M., Tran, K., Radhakrishnan, D., Robinson, A.S., Ogunnaike, B.A., 2014b. Controllability analysis of protein glycosylation in CHO cells. PLOS ONE 9, e87973. doi:10.1371/journal.pone.0087973.

[42] Jimenez del Val, I., Nagy, J.M., Kontoravdi, C., 2011. A dynamic mathematical model for monoclonal antibody N-linked glycosylation and nucleotide sugar donor transport within a maturing Golgi apparatus. Biotechnology Progress 27, 1730–1743. doi:10.1002/btpr.688.

[43] Van Rossum, G., Drake, F.L., 2009. Python 3 Reference Manual. CreateSpace.

[44] Vassiliadis, V.S., Sargent, R.W.H., Pantelides, C.C., 1994. Solution of a class of multistage dynamic optimization problems. 2. Problems with path constraints. Industrial & Engineering Chemistry Research 33, 2123–2133. doi:10.1021/ie00033a015.

[45] Villiger, T.K., Scibona, E., Stettler, M., Broly, H., Morbidelli, M., Soos, M., 2016. Controlling the time evolution of mAb N-linked glycosylation—Part II: Model-based predictions. Biotechnology Progress 32, 1135–1148. doi:10.1002/btpr.2315.

[46] Walsh, G., Walsh, E., 2022. Biopharmaceutical benchmarks 2022. Nature Biotechnology 40, 1722–1760. doi:10.1038/s41587-022-01582-x.

[47] Wong, D.C.F., Wong, K.T.K., Goh, L.T., Heng, C.K., Yap, M.G.S., 2005. Impact of dynamic online fed-batch strategies on metabolism, productivity and N-glycosylation quality in CHO cell cultures. Biotechnology and Bioengineering 89, 164–177. doi:10.1002/bit.20317.

[48] Zupke, C., Brady, L.J., Slade, P.G., Clark, P., Caspary, R.G., Livingston, B., Taylor, L., Bigham, K., Morris, A.E., Bailey, R.W., 2015. Real-time product attribute control to manufacture antibodies with defined N-linked glycan levels. Biotechnology Progress 31, 1433–1441. doi:10.1002/btpr.2136.

**Nonlinear Model Predictive Control of the Glycosylation Process of Monoclonal Antibody in CHO Cell Culture**

Yingjie Ma, Jing Guo, Alexis B. Dubs, and Richard D. Braatz

Massachusetts Institute of Technology, Cambridge, MA, 02139, USA

Table S1. The nucleotide concentration in the Golgi apparatus.

| Substrates | UDP | GDP | CMP |
|---|---|---|---|
| Concentration (mM) | 1.942 | 0.496 | 0.248 |

UDP: Uridine diphosphate, GDP: Guanosine diphosphate, CMP: Cytidine monophosphate.

Table S2. Ground-truth and initial cell culture model parameter values used in ANMPC.

| Variable | Ground truth | Init_0 | Init_1 | Init_2 | Init_3 | Init_4 |
|---|---|---|---|---|---|---|
| $\mu_{\mathrm{max}}$ | −2.73337 | −2.55989 | −2.80766 | −2.79421 | −2.85756 | −2.61458 |
| $\mu_{\mathrm{death,max}}$ | −4.19971 | −5.46269 | −3.05622 | −4.73958 | −4.0814 | −4.3067 |
| $K_{\mathrm{Glc}}$ | 2.64191 | 6.39693 | −1.96326 | −1.96326 | −1.96326 | 6.39693 |
| $K_{\mathrm{Gal}}$ | 2.903069 | 2.559989 | 2.856218 | 2.616629 | 2.919065 | 3.123035 |
| $K_{\mathrm{Urd}}$ | 1.94591 | 1.794221 | 2.103352 | 2.024514 | 2.011801 | 2.020159 |
| $K_{\mathrm{Asn}}$ | 0.963174 | 0.735057 | 1.213042 | 1.00316 | 1.046558 | 0.79888 |
| $\mathrm{KI}_{\mathrm{Amm}}$ | 1.153732 | 0.988067 | 1.253167 | 1.185626 | 1.228392 | 1.050346 |
| $\mathrm{KI}_{\mathrm{Urd}}$ | 3.715765 | 3.834552 | 3.546063 | 3.671895 | 3.732215 | 3.755688 |
| $\mathrm{Kd}_{\mathrm{Amm}}$ | 2.65886 | 1.624842 | 3.543224 | 2.105782 | 3.017001 | 2.680639 |
| $\mathrm{Kd}_{\mathrm{Urd}}$ | 3.327192 | 1.653394 | 4.983583 | 2.841407 | 3.344137 | 3.205529 |
| $Y_{\mathrm{mAb}_X}$ | 1.22083 | 1.246199 | 1.197321 | 1.190082 | 1.22175 | 1.238312 |
| $m_{\mathrm{mAb}}$ | −0.8916 | −0.89174 | −0.87547 | −0.86846 | −0.89231 | −0.88482 |
| $Y_{X_{\mathrm{Glc}}}$ | 20.73322 | 20.69026 | 20.79923 | 20.75349 | 20.71485 | 20.72578 |
| $Y_{X_{\mathrm{Gln}}}$ | 22.25798 | 22.23582 | 22.29092 | 22.27095 | 22.27096 | 22.22578 |
| $Y_{X_{\mathrm{Lac}}}$ | 17.81554 | 17.75285 | 17.80856 | 17.83454 | 17.86724 | 17.80636 |
| $Y_{X_{\mathrm{Amm}}}$ | 21.58193 | 21.56213 | 21.58788 | 21.58079 | 21.57993 | 21.56905 |
| $Y_{X_{\mathrm{Glu}}}$ | 23.40429 | 23.38461 | 23.40029 | 23.40512 | 23.41202 | 23.39805 |
| $Y_{X_{\mathrm{Gal}}}$ | 18.74276 | 19.01825 | 18.75052 | 18.98836 | 18.72869 | 18.60622 |
| $Y_{X_{\mathrm{Urd}}}$ | 21.1995 | 21.28306 | 21.10624 | 21.21614 | 21.22185 | 21.24008 |
| $Y_{X_{\mathrm{Asn}}}$ | 20.4593 | 20.43981 | 20.46797 | 20.45803 | 20.45603 | 20.44601 |
| $m_{\mathrm{Glc}}$ | −24.0959 | −24.1039 | −24.0552 | −24.0824 | −24.1036 | −24.0954 |
| $m_{\mathrm{Lac}}$ | −22.3999 | −23.6157 | −22.7132 | −22.2846 | −21.33 | −22.1211 |
| $K_{C_{\mathrm{Gal}}}$ | 1.66203 | 1.637797 | 1.608433 | 1.646327 | 1.644553 | 1.741433 |
| $f_{\mathrm{Gal}}$ | −1.04982 | −0.95896 | −1.03361 | −1.0096 | −1.00322 | −1.07413 |
| $\mathrm{Lac}_{\mathrm{max1}}$ | 3.054001 | 3.041958 | 3.051802 | 3.062378 | 3.069329 | 3.047592 |
| $\mathrm{Lac}_{\mathrm{max2}}$ | 2.772589 | 2.183389 | 2.626625 | 2.808422 | 3.273907 | 2.950022 |
| $Y_{\mathrm{Gln/Amm}}$ | −2.30259 | −2.77241 | −0.5385 | −1.60263 | −1.43506 | −3.31356 |
| $Y_{\mathrm{Asn/Asp}}$ | −2.30259 | −2.30262 | −2.30257 | −2.3026 | −2.30258 | −2.3026 |
| $Y_{\mathrm{Amm/Urd}}$ | 0.693147 | 0.666049 | 0.744781 | 0.766537 | 0.751395 | 0.770753 |

Table S3. The initial NSD model parameter values used in ANMPC.

| Variable | Ground truth | Init_0 | Init_1 | Init_2 | Init_3 | Init_4 |
|---|---|---|---|---|---|---|
| $V_{\mathrm{max,1}}$ | 0.07851 | 2.221375 | −2.3838 | −2.20738 | −2.3838 | 2.221375 |
| $V_{\mathrm{max,2}}$ | −4.16763 | −4.92916 | −3.54846 | −4.44524 | −4.08572 | −4.24891 |
| $V_{\mathrm{max,2b}}$ | 4.097821 | 6.388225 | 1.783055 | 1.783055 | 1.783055 | 6.388225 |
| $V_{\mathrm{max,3}}$ | −2.937 | −3.92764 | −3.11023 | −3.74464 | −2.90283 | −2.3039 |
| $V_{\mathrm{max,4}}$ | −3.49024 | −3.52893 | −3.42217 | −3.43708 | −3.47991 | −3.39183 |
| $V_{\mathrm{max,5}}$ | −9.20488 | −9.23247 | −9.20829 | −9.23741 | −9.2117 | −9.18441 |
| $V_{\mathrm{max,6}}$ | 1.606501 | −0.66748 | 3.937691 | −0.66748 | −0.66748 | −0.66748 |
| $V_{\mathrm{max,7}}$ | 1.572593 | 2.523005 | −0.77653 | 3.828641 | 3.828641 | 3.828641 |
| $V_{\mathrm{max,1_{Urd}}}$ | −1.92801 | −1.92909 | −1.97474 | −1.95722 | −1.83945 | −1.91091 |
| $V_{\mathrm{max,2_{Urd}}}$ | −2.86122 | −2.90539 | −2.85582 | −2.74666 | −2.8638 | −2.82725 |
| $V_{\mathrm{max,4_{Urd}}}$ | −4.38414 | −4.41833 | −4.24971 | −4.42033 | −4.56943 | −4.41704 |
| $V_{\mathrm{max,6_{Urd}}}$ | 1.69449 | 3.35425 | 3.978316 | 3.978316 | 2.500113 | 3.978316 |
| $V_{\mathrm{max,7_{sink}}}$ | 2.349945 | 4.695011 | 4.695011 | 4.695011 | 0.089841 | 0.089841 |
| $V_{\mathrm{max,1_{sink}}}$ | 3.432078 | 0.935701 | 5.540871 | 5.540871 | 5.540871 | 0.935701 |
| $V_{\max,6_{\text{sink}}}$ | 2.0374 | 3.687171 | 4.290459 | 4.290459 | 2.84107 | 4.290459 |
| $V_{\max,6_{\text{Gal}}}$ | 3.638536 | 4.704453 | 5.661533 | 5.905233 | 4.216692 | 6.02519 |
| $f_{\text{Gln}}$ | −3.79203 | −6.21461 | −1.60944 | −1.60944 | −1.60944 | −6.21461 |
| $K_{\text{M1}_{\text{Gln}}}$ | −0.9875 | −3.17009 | 1.435085 | 1.435085 | 1.435085 | −3.17009 |
| $K_{\text{M1}_{\text{sink}}}$ | −3.21452 | −3.08029 | −3.08511 | −3.1052 | −2.9484 | −3.33346 |
| $\text{KI}_{1_{\text{sink}}}$ | −8.82695 | −6.71713 | −11.3223 | −11.3223 | −11.3223 | −6.71713 |
| $K_{\text{M2}_{\text{Glc}}}$ | 4.394494 | 2.055661 | 6.660831 | 2.98338 | 5.371626 | 3.593824 |
| $K_{\text{M2b}_{\text{UDPGal}}}$ | −3.70909 | −1.39433 | −5.9995 | −5.9995 | −5.9995 | −1.39433 |
| $\text{KI}_{2\text{A}}$ | −13.8113 | −16.0693 | −11.4641 | −11.4641 | −11.4641 | −16.0693 |
| $\text{KI}_{2\text{B}}$ | 4.522982 | 4.775267 | 4.136864 | 4.457755 | 4.455283 | 4.759863 |
| $\text{KI}_{2\text{C}}$ | −4.31993 | −6.62258 | −2.01741 | −2.01741 | −2.01741 | −6.62258 |
| $\text{KI}_{2\text{D}}$ | −12.7805 | −15.1398 | −10.5346 | −10.5346 | −10.5346 | −15.1398 |
| $K_{\text{M3}_{\text{Glc}}}$ | 3.938518 | 5.046892 | 4.912108 | 4.837086 | 3.755154 | 3.480013 |
| $K_{\text{M4}_{\text{UDPGlcNAc}}}$ | 0.969506 | 0.917998 | 1.05351 | 1.021446 | 0.987863 | 1.081214 |
| $K_{\text{M5}_{\text{UDPGlcNAc}}}$ | −3.59839 | −3.92842 | −3.536 | −3.52125 | −3.57602 | −3.16815 |
| $\text{KI}_5$ | 6.907754 | 9.21034 | 4.60517 | 4.60517 | 4.60517 | 9.21034 |
| $K_{\text{M6}_{\text{UDPGlc}}}$ | −4.10657 | −1.83258 | −1.83258 | −1.83258 | −1.83258 | −6.43775 |
| $\text{KI}_{6\text{A}}$ | −16.0402 | −13.7202 | −18.3254 | −13.7202 | −13.7202 | −13.7202 |
| $\text{KI}_{6\text{B}}$ | 1.519513 | 2.110467 | 0.642307 | 1.384974 | 1.361761 | 2.067355 |
| $\text{KI}_{6\text{C}}$ | −12.2498 | −14.5412 | −14.5412 | −14.5412 | −14.5412 | −14.5412 |
| $K_{\text{M6}_{\text{sink}}}$ | −2.05652 | −3.95876 | −4.34281 | −4.34281 | −2.85926 | −4.34281 |
| $\text{KI}_{6_{\text{sink}}}$ | −1.12932 | −0.95854 | −0.92372 | −0.88972 | −1.09845 | −1.09435 |
| $K_{\text{M7}_{\text{GDPMan}}}$ | −0.05659 | 0.443242 | 0.478822 | 2.292535 | 2.292535 | −0.14687 |
| $\text{KI}_7$ | −4.10082 | −6.41306 | −1.80789 | −3.59805 | −6.41306 | −6.41306 |
| $K_{\text{M7}_{\text{sink}}}$ | 2.226197 | 4.486387 | 4.486387 | 4.486387 | −0.11878 | −0.11878 |
| $K_{M1_{\text{Urd}}}$ | 1.889027 | 1.761356 | 1.881457 | 1.926024 | 2.064203 | 1.914329 |
| $K_{M2_{\text{Urd}}}$ | 2.396853 | 2.264824 | 2.457587 | 2.80713 | 2.527026 | 2.424083 |
| $K_{M4_{\text{Urd}}}$ | 1.872129 | 1.865431 | 1.994391 | 1.842122 | 1.592921 | 1.897788 |
| $K_{M6_{\text{Urd}}}$ | −0.84622 | −2.69789 | −3.96703 | −4.12308 | −1.69607 | −4.35423 |
| $K_{M6_{\mathrm{Gal}}}$ | −0.44656 | −2.78634 | −4.09344 | −4.70629 | −0.93837 | −4.81493 |
| $\mathrm{KI}_{6\mathrm{D}}$ | −4.67068 | −4.45774 | −4.593 | −5.10185 | −4.04721 | −4.96951 |
| $\mathrm{KI}_{6\mathrm{E}}$ | 4.601409 | 5.180291 | 3.644233 | 4.405577 | 4.449866 | 5.122605 |
| $\mathrm{KI}_{6\mathrm{F}}$ | −6.99762 | −9.30355 | −9.30355 | −9.30355 | −7.97132 | −9.30355 |
| $K_{\mathrm{TP,UDP-Glc}}$ | −0.0266 | 0.079036 | −0.0399 | 0.032919 | 0.109992 | −0.08821 |
| $K_{\mathrm{TP,UDP-GlcNAc}}$ | 1.620836 | 3.726258 | 0.480274 | 0.723842 | −0.6832 | 3.921973 |
| $K_{\mathrm{TP,UDP-Gal}}$ | 1.93832 | 0.326087 | 2.172851 | 1.172591 | 1.40261 | 0.978521 |
| $K_{\mathrm{TP,UDP-GalNAc}}$ | 2.414956 | 2.359934 | 2.382138 | 2.287052 | 2.419422 | 2.505055 |
| $K_{\mathrm{TP,GDP-Man}}$ | −2.06746 | −0.2379 | −1.94578 | −1.09393 | −1.81769 | −2.66049 |
| $K_{\mathrm{TP,GDP-Fuc}}$ | −2.30651 | −4.60517 | 0 | 0 | −4.60517 | −4.60517 |
| $K_{\mathrm{TP,CMP-Neu5Ac}}$ | 6.217368 | 6.105463 | 6.257868 | 6.379847 | 6.20563 | 6.534773 |

Table S4. The initial Golgi model parameter values used in ANMPC.

| Variable | Ground truth | Init_0 | Init_1 | Init_2 | Init_3 | Init_4 |
|---|---|---|---|---|---|---|
| $K_{d,\mathrm{OS1}}^{\mathrm{ManI}}$ | 4.075253 | 6.405228 | 1.800058 | 1.800058 | 1.800058 | 6.405228 |
| $K_{d,\mathrm{OS2}}^{\mathrm{Man1}}$ | 4.638963 | 2.397895 | 7.003065 | 2.397895 | 2.397895 | 7.003065 |
| $K_{d,\mathrm{OS3}}^{\mathrm{Man1}}$ | 3.408208 | 1.12493 | 5.7301 | 1.12493 | 5.7301 | 1.12493 |
| $K_{d,\mathrm{OS4}}^{\mathrm{Man1}}$ | 4.251972 | 6.608001 | 2.00283 | 6.608001 | 6.608001 | 2.00283 |
| $K_{d,\mathrm{OS6}}^{\mathrm{Man2}}$ | 5.109051 | 2.995732 | 7.600902 | 2.995732 | 7.600902 | 2.995732 |
| $K_{d,\mathrm{OS7}}^{\mathrm{Man2}}$ | 4.513693 | 6.907755 | 2.302585 | 6.907755 | 2.302585 | 6.907755 |
| $K_{d,\mathrm{OS5}}^{\mathrm{GnT1}}$ | 4.384545 | 6.870053 | 6.870053 | 2.264883 | 2.264883 | 2.264883 |
| $K_{d,\mathrm{OS9}}^{\mathrm{GnT2}}$ | 3.958708 | 6.877296 | 2.272126 | 6.877296 | 6.877296 | 2.272126 |
| $K_{d,\mathrm{OS6}}^{\mathrm{FucT}}$ | 3.942005 | 6.663133 | 2.057963 | 4.292903 | 2.667411 | 4.85936 |
| $K_{d,\mathrm{OS9}}^{\mathrm{GalT}}$ | 8.498685 | 8.483642 | 8.483617 | 8.501296 | 8.508624 | 8.499496 |
| $K_{d,\mathrm{OS12}}^{\mathrm{GalT}}$ | 7.462721 | 7.469098 | 7.454398 | 7.461404 | 7.459033 | 7.467818 |
| $K_{d,\mathrm{OS13}}^{\mathrm{SiaT}}$ | 10.57627 | 12.85055 | 8.245384 | 8.245384 | 8.245384 | 12.85055 |
| $K_{d,\mathrm{Mn}}^{\mathrm{GnT1}}$ | −5.20985 | −7.51106 | −7.51106 | −2.90589 | −2.90589 | −6.98216 |
| $K_{d,\mathrm{Mn}}^{\mathrm{GnT2}}$ | −5.21677 | −5.54648 | −5.06692 | −5.05458 | −4.93477 | −5.53594 |
| $K_{d,\mathrm{Mn}}^{\mathrm{GalT}}$ | −2.13271 | −5.5675 | −5.5675 | −5.5675 | −0.96233 | −5.5675 |
| $K_{d,\text{UDP-GlcNAc}}^{\text{GnT1}}$ | 5.113426 | 2.833213 | 2.833213 | 7.438384 | 7.438384 | 2.833213 |
| $K_{d,\text{UDP-GlcNAc}}^{\text{GnT2}}$ | 6.354573 | 4.564348 | 9.169518 | 9.169518 | 9.169518 | 4.564348 |
| $K_{d,\text{GDP-Fuc}}^{\text{FucT}}$ | 3.563574 | 1.526056 | 6.131226 | 2.803191 | 5.444407 | 3.083855 |
| $K_{d,\text{UDP-Gal}}^{\text{GalT}}$ | 5.616017 | 6.476972 | 6.476972 | 6.476972 | 2.943144 | 6.476972 |
| $K_{d,\text{CMP-Neu5Ac}}^{\text{SiaT}}$ | 3.91967 | 1.609438 | 6.214608 | 6.214608 | 6.214608 | 1.609438 |

Table S5. ANMPC results with and without a bleed stream (Init_1; partial measurements; $\Delta t_{\text{prep}} = 4$ h).

| | Viability at harvest (%) | GI (mg L$^{-1}$) | Constraint status |
|---|---|---|---|
| No bleed | 60.75 | 395.25 | Feasible |
| With bleed | 65.90 | 368.28 | Feasible |

Table S6. ANMPC results with and without a bleed stream (Init_2; partial measurements; $\Delta t_{\text{prep}} = 4$ h).

| | Viability at harvest (%) | GI (mg L$^{-1}$) | Constraint status |
|---|---|---|---|
| No bleed | 57.21 | 415.79 | Infeasible |
| With bleed | 72.07 | 279.19 | Feasible |

Table S7. ANMPC results with and without a bleed stream (Init_3; partial measurements; $\Delta t_{\text{prep}} = 4$ h).

| | Viability at harvest (%) | GI (mg L$^{-1}$) | Constraint status |
|---|---|---|---|
| No bleed | 56.87 | 417.93 | Infeasible |
| With bleed | 65.27 | 363.95 | Feasible |

Table S8. ANMPC results with and without a bleed stream (Init_4; partial measurements; $\Delta t_{\text{prep}} = 4$ h).

| | Viability at harvest (%) | GI (mg L$^{-1}$) | Constraint status |
|---|---|---|---|
| No bleed | 56.37 | 415.74 | Infeasible |
| With bleed | 64.80 | 374.16 | Feasible |